\documentclass[a4paper,oneside,english,reqno,12pt]{amsart}

\usepackage[utf8]{inputenx}
\usepackage[russian,english]{babel}

\usepackage{amssymb,amsfonts,amsthm,amsmath,mathrsfs,amscd,comment,latexsym,bbold,dsfont}
\usepackage{amsthm, mathtools, stmaryrd}
\usepackage{textcomp, url, enumerate, latexsym, graphicx, varioref, hyperref, multirow, layout, siunitx, booktabs}

\usepackage{blindtext}

\usepackage{stackengine}

\usepackage[matrix,arrow,curve]{xy}

\usepackage{enumitem}

\usepackage[top=2cm,left=2cm,bottom=2.5cm,right=2cm]{geometry}

\RequirePackage[backend    = biber,
                sortcites  = true,
                giveninits = true,
                doi        = false,
                isbn       = false,
                url        = false,
                maxnames   = 6,
                style      = numeric]{biblatex}
\DeclareFieldFormat*{title}{#1}
\usepackage{cleveref}

\numberwithin{equation}{section}

\theoremstyle{plain}
\newtheorem{theorem}{Theorem}
\newtheorem{corollarytheorem}[theorem]{Theorem}
\newtheorem*{theorem*}{Theorem}

\newtheorem{corollary}[theorem]{Corollary}
\newtheorem*{corollary*}{Corollary}

\newtheorem{lemma}{Lemma}
\newtheorem{lemmacorollary}[lemma]{Corollary}

\theoremstyle{definition}
\newtheorem{definition}{Definition}

\newtheorem{notation}{Notation}
\newtheorem{definition-notation}{Definition-Notation}

\theoremstyle{remark}
\newtheorem{remark}{Remark}

\theoremstyle{plain}
\newtheorem*{nonum-theorem}{Theorem}

\theoremstyle{plain} 
\newcommand{\thistheoremname}{}
\newtheorem*{genericthm}{\thistheoremname}

\newtheoremstyle{named}{}{}{\itshape}{}{\bfseries}{.}{.5em}{\thmnote{#3's }#1}

\numberwithin{theorem}{section}
\numberwithin{lemma}{section}
\numberwithin{proposition}{section}
\numberwithin{definition}{section}

\newcommand{\id}{\mathrm{id}}

\newcommand\barbelow[1]{\stackunder[1.2pt]{$#1$}{\rule{1ex}{.15ex}}}
\newcommand{\overbar}[1]{\mkern 1.5mu\overline{\mkern-1.5mu#1\mkern-1.5mu}\mkern 1.5mu}

\newcommand{\Sch}{\mathrm{Sch}}
\newcommand{\Sm}{\mathrm{Sm}}

\newcommand{\calL}{\mathcal{L}}
\newcommand{\calS}{\mathcal{S}}

\newcommand{\climmersion}{c}

\newcommand{\ovX}{\overline{X}}

\newcommand{\Monrho}{{\rho^\Mon}}
\newcommand{\Monrhoprime}{\rhoSetprime^\Mon}
\newcommand{\rhoSet}{\rho}
\newcommand{\rhoSetP}{\rho_P}
\newcommand{\element}{\mathrm{el}}
\newcommand{\rhoelement}{\rho^\element}
\newcommand{\rhoMonElement}{\rho^{\Mon\text{-}\element}}
\newcommand{\rhoMonSubset}{\rho^\mathrm{subset}}
\newcommand{\rhoElement}{\rhoelement}
\newcommand{\rhoelementprime}{{\barbelow{\rho}^{\element}}}
\newcommand{\rhoElementprime}{\rhoelementprime}
\newcommand{\rhoSetprime}{\barbelow{\rho}}

\newcommand{\oX}{X}
\newcommand{\schX}{\mathsf{X}}
\newcommand{\schx}{\mathsf{x}}
\newcommand{\ox}{x}
\newcommand{\oXprime}{\barbelow{\oX}}
\newcommand{\schXprime}{\barbelow{\mathsf{X}}}

\newcommand{\aroX}{\vec{X}}

\newcommand{\gr}{\partial} 
\newcommand{\supp}{\varsigma} 

\newcommand{\ras}{r}

\newcommand{\grprime}{\barbelow{\gr}}
\newcommand{\suppprime}{\barbelow{\supp}}
\newcommand{\catO}{O}
\newcommand{\precatO}{O^c}
\newcommand{\catAO}{\vec{\mathcal{A}}}
\newcommand{\catAOu}{\mathcal{A}}
\newcommand{\precatAOu}{A}
\newcommand{\precatAO}{\vec{A}}
\newcommand{\precatOB}{\precatO_B}

\newcommand{\catOB}{\catO_B}
\newcommand{\catOBregstratasgrsupp}{\catOB}

\newcommand{\precatAOB}{\precatAO_B}
\newcommand{\precatAOuB}{\precatAOu_B}
\newcommand{\catAOB}{\catAO_B}
\newcommand{\catAOuB}{\catAOu_B}
\newcommand{\FunctorSch}{\FschX}
\newcommand{\funcXcomplinfX}{{\Sch^V}}
\newcommand{\funcXcomplinfXquotientXcomplinfXunionSup}{{\Sch^Q}}

\newcommand{\compinttimmersion}{\xi}

\newcommand{\limBlras}{\mathrm{limBl}_\ras}

\newcommand{\limBlbeta}{\mathrm{limBl}_\beta}
\newcommand{\Blras}{{\Bl_\ras}}
\newcommand{\Blrasequiv}{r}
\newcommand{\Blrho}{{\Bl_{\rho}}}
\newcommand{\Blbeta}{{\Bl_{\beta}}}

\newcommand{\LBlbeta}{L_{\beta}} 
\newcommand{\SpccompintABlbeta}{\Spc_{\compintA,\beta}}
\newcommand{\uLLP}{\mathrm{uLLP}} 
\newcommand{\uRLP}{\mathrm{uRLP}}
\newcommand{\Et}{\mathrm{Et}}
\newcommand{\EtInfNeigh}{\mathrm{InfEtNeigh}}

\newcommand{\Bl}{{Bl}}
\newcommand{\preimagerasrho}{\preimagerho}

\newcommand{\preimagegrsup}{$\eta$} 
\newcommand{\rhoirr}{$\nu$}

\newcommand{\infinitesimaletaleneighbourhood}{infinitesimal \'etale neighbourhood}
\newcommand{\intA}{{\A^1}}
\newcommand{\LintA}{L_{\intA}}
\newcommand{\compintA}{\underline{\overline{\A}}}
\newcommand{\compintAone}{\underline{\overline{\A^1}}}
\newcommand{\DeltacompintA}[1]{\Delta^{#1}_{\compintA}}
\newcommand{\DeltacompintAall}[1]{\Delta^{#1}_{\compintA,\all}}
\newcommand{\DeltaonecompintA}{\Delta^1_{\compintA}}
\newcommand{\DeltaonecompintAall}{\Delta^1_{\compintA,\all}}
\newcommand{\compintAall}{\compintA^1_{\all}}
\newcommand{\DeltalcompintA}{\Delta^{\times l}_{\compintA}}
\newcommand{\DeltazerocompintA}{\Delta^{0}_{\compintA}}
\newcommand{\DeltabonecompintA}{\Delta^{\bullet\leq 1}_{\compintA}}
\newcommand{\DeltapcompintA}{\Deltatimesbullet_{\compintA}}
\newcommand{\all}{\mathrm{all}}
\newcommand{\A}{\mathbb{A}}
\newcommand{\PP}{\mathbb{P}}\newcommand{\rhoblowup}{$\rho$-blow-up}
\newcommand{\betablowup}{$\beta$-blow-up}
\newcommand{\elementaryrhoblowup}{elementary $\rho$-blow-up}
\newcommand{\elementarybetablowup}{elementary $\beta$-blow-up}
\newcommand{\rhoblowups}{$\rho$-blow-ups}
\newcommand{\betablowups}{$\beta$-blow-ups}

\newcommand{\LcompintABlbeta}{L_{\compintA,\beta}}

\newcommand{\LBlbetacompintA}{\LcompintABlbeta}

\newcommand{\geomrealKcompintA}{{|K|_{\compintA}}}
\newcommand{\geomrealcompactintAK}{{ |K|_{\compintA,\all} }}
\newcommand{\geomrealcompactintASn}{{ |S^n|_{\compintA,\all} }}

\newcommand{\geomrealallrascompactintAK}{\geomrealcompactintAK}

\newcommand{\ABHBC}{\phi}
\newcommand{\ABHBCprop}{\ABHBC}

\newcommand{\RLPwrt}{satisfies RLP w.r.t. }

\newcommand{\compintABlbeta}{{\compintA,\Blbeta}}
\newcommand{\trivcofHBC}{(\ABHBC,\compintA,\beta)} 
\newcommand{\compintAoneBlbeta}{(\compintA,\beta)} 

\newcommand{\Compactified}{\mathrm{comp}}
\newcommand{\SchCompactified}{\Sch^\Compactified}
\newcommand{\SchComp}{\SchCompactified}
\newcommand{\limSchComp}{\mathrm{lim}\Sch^\Compactified}
\newcommand{\SchCompactifiedRegular}{\Sch^\Compactified_\mathrm{reg}}
\newcommand{\FunctorSchCompactified}{\SchCompactified}
\newcommand{\funcXcomplinfXquotientXcomplinfXunionSupCompactified}{Q^\Compactified}

\newcommand{\Spcpointed}{\Spc^\bullet}

\newcommand{\topology}{\tau}
\newcommand{\Nis}{\mathrm{Nis}}
\newcommand{\Zar}{\mathrm{Zar}}
\newcommand{\Spc}{\mathrm{Spc}}
\newcommand{\Spt}{\mathrm{Spt}}
\newcommand{\SpcpresheafonSchB}{W}
\newcommand{\Hmot}{\mathrm{H}}
\newcommand{\SH}{\mathrm{SH}} %
\newcommand{\SHSone}{\SH^{S^1}}
\newcommand{\Hpointed}{\Hmot^\bullet}
\newcommand{\HpointedAtopology}{\Hpointed_{\A^1,\topology}}
\newcommand{\SpcAtopology}{\Spc_{\A^1,\topology}}

\newcommand{\HpointedANis}{\Hpointed_{\A^1,\Nis}}

\newcommand{\SHAZar}{\SH_{\A^1,\Zar}}

\newcommand{\motSpt}{S}
\newcommand{\motSpc}{S}

\newcommand{\motSpectrum}{S}
\newcommand{\PresheafmotSpectrum}{\mathrm{S}}
\newcommand{\Map}{\mathrm{Map}}
\newcommand{\MapSet}{\mathrm{Map}}
\newcommand{\MapSpc}{\mathrm{Map}^\Spc}
\newcommand{\MapSpt}{\mathrm{Map}^\Spt}

\newcommand{\HonSpt}{H}
\newcommand{\HKonSpt}{H^{K}}
\newcommand{\HKonSpttau}{\HKonSpt_\tau}

\newcommand{\HSnonSpt}{H^{S^n}}

\newcommand{\HSnonSptNis}{\HSnonSpt_\Nis}
\newcommand{\HonSpc}{H}
\newcommand{\HKonSpc}{H^{K}}
\newcommand{\HKonSpctau}{\HKonSpc_\tau}

\newcommand{\HSnonSpc}{H^{S^n}}

\newcommand{\HSnonSpcNis}{\HSnonSpc_\Nis}
\newcommand{\FunctorSpcSmBtoSpcSchB}{E^\Spc}
\newcommand{\FunctorHpointedAtopologySmBtoHpointedAtopologySchB}{E^{\Hmot}}

\newcommand{\eqC}{ec}

\newcommand{\ovoX}{\overbar{X}}

\DeclareMathOperator{\codim}{codim}

\newcommand{\Mon}{\mathrm{Mon}}

\newcommand{\Func}{\mathrm{Func}}
\newcommand{\op}{\mathrm{op}}

\newcommand{\diagramBlrhoeq}{D_{\rho\text{-bl-eq}}}
\newcommand{\diagramBlbetaeq}{D_{\beta\text{-bl-eq}}}

\newcommand{\cubicalset}{cubical set}

\newcommand{\Set}{\mathrm{Set}}

\newcommand{\Notationref}[1]{Notation \ref{#1}}

\newcommand{\Zprime}{Z^\prime}

\newcommand{\subsetrhoSet}{\rho^{\mathrm{Subset}}}
\newcommand{\BloX}{\tilde{\oX}}
\newcommand{\sinLcapZsequalsZeqmpty}{{sL}^{\emptyset}}
\newcommand{\sinLcapZsequalsZempty}{\sinLcapZsequalsZeqmpty}
\newcommand{\sinLcapZsequalsZdelta}{sL^{\mathrm{div}}}
\newcommand{\sinLcapZssubsetcupZs}{sL^{\mathrm{\subset}}}
\newcommand{\sinLonXcapZssubsetcupZs}{\sinLcapZssubsetcupZs}
\newcommand{\sLonXKZssubsetcupZs}{\sinLcapZssubsetcupZs}

\newcommand{\supextensionmorphismHANisB}{E_\Nis}
\newcommand{\supextensionmorphismtU}{E}
\newcommand{\supextensionmorphismtUtilde}{\tilde{E}}
\newcommand{\objectUrho}{\mathcal{U}}
\newcommand{\objectUrhotilde}{\tilde{\objectUrho}} 
\newcommand{\supextensionmorphismHAtopologyB}{E_\topology}

\newcommand{\timmersionUrho}{\tilde{\climmersion}}

\newcommand{\regularstratarhoprecatAOuB}{\precatAOu_{\mathrm{rs},B}}
\newcommand{\regularstratarhoprecatAOB}{\precatAO_{\mathrm{rs},B}}
\newcommand{\regularstratarhocatAOB}{\catAO_{\mathrm{rs},B}}

\newcommand{\indexSubsets}{\alpha}
\newcommand{\indexSubsetsempty}{\epsilon} 
\newcommand{\indexSubsetssrho}{\gamma}
\newcommand{\SetofindexesSubsets}{A}

\newcommand{\rhomap}{{\rho^\mathrm{Map}}}

\newcommand{\bl}{bl}

\newcommand{\Deltableqone}{\Delta^{\bullet\leq 1}}
\newcommand{\Deltatimesbullet}{\Delta^{\times\bullet}}
\newcommand{\Deltatimesb}{\Deltatimesbullet}

\newcommand{\Gr}{\mathrm{Gr}}

\newcommand{\LAZar}{L_{\A^1,\Zar}}
\newcommand{\LANis}{L_{\A^1,\Nis}}

\newcommand{\rhoPic}{$\mathscr{Pic}$}

\newcommand{\rhopreimage}{\preimagerho}

\newcommand{\good}{good}

\newcommand{\rhoSubset}{\subsetrhoSet}

\newcommand{\classa}{\mathbf{a}}
\newcommand{\rasconstmorphismsclassa}{\classa_r}
\newcommand{\rasmorphismsclassa}{\classa^r}
\newcommand{\constroveroX}{\breve{\oX}}
\newcommand{\constroXprime}{\barbelow{\breve{\oX}}}
\newcommand{\constroveroXprime}{\constroXprime}
\newcommand{\preimagerho}{$\mu$}
\newcommand{\rhoimmersion}{$\pi$}
\newcommand{\rhotriv}{$\tau$}
\newcommand{\rliftrho}{r}

\newcommand{\elementinsetsSchrhoBl}{\beta^\mathrm{el}}

\newcommand{\SubsetMonrho}{\beta^\mathrm{el}}
\newcommand{\betaElement}{\beta^\mathrm{el}}

\newcommand{\presetsSchrhoBl}{\beta}
\newcommand{\presetsSchrhoBlprime}{\barbelow{\presetsSchrhoBl}}
\newcommand{\presetsSchrhoBlSubset}{\presetsSchrhoBl^\mathrm{Subset}}
\newcommand{\cupofpresetsSchrhoBlSubset}{\cup(\presetsSchrhoBlSubset)}
\newcommand{\cupofpresetsSchrhoBlSubsetSubset}{c}

\newcommand{\setsSchrhoBl}{\presetsSchrhoBlcup^{\Mon}}
\newcommand{\setsSchrhoBlprime}{\presetsSchrhoBlprime^{\Mon}}
\newcommand{\presetsSchrhoBlcupprime}{\presetsSchrhoBlcup^\prime}

\newcommand{\rhoci}{$\rho$-ci}
\newcommand{\FschX}{\Sch^U_B}
\newcommand{\precatAOuBtoSch}{\FschX}
\newcommand{\eqXinovX}{e}

\newcommand{\schovX}{\overline{\schX}}
\newcommand{\ovschX}{\schovX}
\newcommand{\byequations}{equational}
\newcommand{\eX}{V}
\newcommand{\ecompltoX}{C}
\newcommand{\etalebyequations}{\EtEquational}
\newcommand{\EtEq}{\mathrm{EtEq}}

\newcommand{\EtEquational}{\'etale-equational}
\newcommand{\CapitalEtEquational}{\'Etale-equational}
\newcommand{\etalemorphism}{e}
\newcommand{\EtEqIRRho}{\mathrm{EtEq}_{\rho}} 
\newcommand{\productprojectivespaceP}{$\mathbb{P}$}
\newcommand{\productprojectivespace}{\productprojectivespaceP-object}
\newcommand{\productprojectivespacemorphism}{\productprojectivespaceP-morphism} 

\newcommand{\limoX}{\tilde{X}}
\newcommand{\rhoblmorphism}{b}
\newcommand{\rhoblmorphismmu}{b}
\newcommand{\limBlrhooX}{B} 
\newcommand{\limElimBlrhooX}{E}
\newcommand{\BlrhooX}{\mathrm{Bl}_\rho(\oX)}

\newcommand{\rasarrows}{\chi} 
\newcommand{\rasarrow}{f}
\newcommand{\rasoX}{{\oX^r}}
\newcommand{\Diagramrasarrows}{D(\oX)}
\newcommand{\LrhoBl}{L_{\Bl_\rho}}
\newcommand{\LbetaBl}{L_{\Bl_\beta}}
\newcommand{\ClassImmersionsHBCforLift}{\mathscr{H}}
\newcommand{\powerSoverL}{\mathfrak{p}}
\newcommand{\powerS}{\mathfrak{S}}
\newcommand{\powerL}{\mathfrak{L}}
\newcommand{\sequation}{e}
\newcommand{\SetofindexSubSets}{A}
\newcommand{\SetofindexesSubsetsrho}{\Gamma}
\newcommand{\arrowequiv}{\shortdownarrow}
\newcommand{\limEBlPobject}{{\mathsf{E}}}
\newcommand{\limBlPobject}{{\mathsf{B}}}

\newcommand{\setsrhoBlOperations}{\cup} 

\newcommand{\Monoid}{Mon}
\newcommand{\Subsets}{\mathrm{Subsets}}

\newcommand{\ssetminus}{-}
\newcommand{\ssetminuszero}[1]{(#1\ssetminus \{0\})}
\newcommand{\setminuszero}[1]{#1_{\setminus \{0\}}}

\newcommand{\setsSchrhoBlElement}{\beta^{\mathrm{el}}}
\newcommand{\setsSchrhoBllimElement}{\beta^{\mathrm{lim}}}

\newcommand{\presetsSchrhoBlf}{\dot{\beta}}
\newcommand{\presetsSchrhoBlfcup}{\dot{\beta}^{\cup}}
\newcommand{\presetsSchrhoBlfcupElement}{\dot{\beta}^{\mathrm{el}}}
\newcommand{\Cov}{\mathrm{Cov}}
\newcommand{\ShiftsColimits}{\mathrm{ShifColim}}

\newcommand{\LZar}{L_{\Zar}}
\newcommand{\LZarintA}{L_{\Zar,\A^1}}

\newcommand{\grcompZar}{\stackrel{\circ}{\mathrm{Zar}}}
\newcommand{\covset}{c}
\newcommand{\betaLiftcalSring}{$\beta$}
\newcommand{\betaLiftcalSringobject}{\betaLiftcalSring-object}
\newcommand{\presetsSchrhoBlcup}{\presetsSchrhoBl}
\newcommand{\coprodpointed}{\coprod}
\newcommand{\Colimits}{\mathrm{Colim}}
\newcommand{\quadforall}{\quad\forall}
\newcommand{\SpherewedgelB}{\PP^l_B/(\PP^l_B-0_B)} 
\newcommand{\rhoSubsetsCentres}{\rho^{SC}}
\newcommand{\rhoSubsetswrtGood}{\rho^{SG}}
\newcommand{\rhoSubsetsGood}{\rhoSubsetswrtGood}
\newcommand{\rhoSubsetCentre}{\rho^{C}} 

\newcommand{\rhoSubsetGood}{\rho^{G}}
\newcommand{\rhoSubsetCentreInductionElement}{\rho^{C}_{f}}
\newcommand{\multa}{a}
\newcommand{\rhoSubsetCentreInductionElementhat}{\hat{\rho}^{C}_{f}}
\newcommand{\multb}{b}
\newcommand{\rhoSubsetsCentreshat}{\hat{\rho}^{SC}}
\newcommand{\seqbl}{bl}
\newcommand{\rhoF}{\rho^{\mathrm{F}}}
\newcommand{\FD}{$\mathscr{F}$}
\newcommand{\exceptional}{e}
\newcommand{\powerexceptional}{p}

\newcommand{\Vect}{\mathrm{Vect}}

\newcommand{\Spec}{\operatorname{Spec}}

\title{Gersten conjecture for K-theory on henselian schemes and $\phi$-motivic localisation}
\author{Andrei~E.~Druzhinin}
\address{
St. Petersburg State University, 14th Line V.O., 29B, Saint Petersburg 199178 Russia,
\\
St. Petersburg Department of V.A.~Steklov Institute of Mathematics of the Russian Academy of Sciences,
191023, 27 Fontanka, St. Petersburg, Russia.
}
\email{andrei.druzh@gmail.com}
\subjclass[MSC2020]{14D99, 14F99}
\keywords{K-theory, Gersten Conjecture, motivic homotopy over DVR}
\begin{document}
\selectlanguage{english}
\maketitle
\begin{abstract}

A key triviality result for 
support extension maps 
for 
motivic $\A^1$-homotopies 
of cellular motivic 
spaces $S$ 
over 
a DVR spectrum $B$
is proven.
Combining with 
earlier known results on 
the K-theory 
Gersten complex 
and 
the K-theory 
motivic spectrum
we achieve a proof of 
the Gersten Conjecture 
on 
essentially smooth local Henselian $B$-schemes.
Additionally, 
we outline 
generalisations 
for 
Cousin complexes associated to 
motivic $\A^1$- and $\square$-homotopies of 
cellular $B$-spectra.

The proof 
is based on two ingredients:
    \begin{itemize}
\item[(1)] 
A new 
``motivic localisation'' over $B$, 
called \emph{$\phi$-motivic}, 
giving
rise to the $\phi$-motivic homotopy category
such that 
the triviality of the support extension maps and the acyclicity of Cousin complexes 
hold
for all 
objects 
$S$, 
not necessarily cellular.

\item[(2)]
An interpretation of 
some 
classes in the motivic $\A^1$-homotopies with support
defined with respect to 
the Morel-Voevodsky motivic homotopy category of smooth $B$-schemes
in terms of 
the construction of $\phi$-motivic homotopy category mentioned 
in Point (1).
\end{itemize}

%
%
%

\end{abstract}\tableofcontents
\section{Introdction}

\subsection{Gersten Conjecture for K-theory} 

\subsubsection{Claim and Results} 
The Gersten Conjecture
\cite{Ge}, \cite[Conjecture 5.10]{zbMATH05778592}
claims the acyclicity of the complex
\begin{equation}\label{eq:KthComplexUn}K_n(U)\to 
\bigoplus_{u\in U^{(0)}}K_{n}(u)\to
\dots \to \bigoplus_{u\in U^{(c)}}K_{n-c}(u)\to \dots\end{equation}
for a local regular scheme $U$ 
\footnote{
claimed after Gersten's result, 
and formulated publically in such generality 
by Quillen.
}
\footnote{
\cite{arXiv:2404.06655} by Kundu 
discusses 
generalisation of the initial morphism injectivity in \eqref{eq:KthComplexUn} for valuation rings.
}
\footnote{
The acyclicity 
holds for $n=0$, $n=1$ and
for $n=2$, 
\cite{arXiv:2202.00896} by Panin.
}%
, where $K_*$ is the
algebraic K-theory
.
The claim is proven 
\begin{itemize}
\item
by Gersten in \cite{Ge},
for 
DVR spectra $U$
with finite residue field%
s, 
\item
extended 
%
to 
finite fields algebraic extensions
in \cite{zbMATH03618257} by Sherman, 
\item 
by Quillen
in \cite{Q}, 
for 
local schemes 
$U=\Spec(O_{X,x})$,
$X\in\Sm_k$,
over a field $k$, 
\item
extended to all local regular equi-characteristic schemes $U$
in \cite{P1} by Panin.
\end{itemize}
Our article proves the claim
\begin{itemize}
\item[$\circ$]
for 
essentially smooth 
local Henselian 
$B$-schemes $U$ over a DVR spectrum $B$.
%
%
\end{itemize}

\begin{remark} 
We remark unpublished preprints by Mochizuki 
referring to the survey \cite{arXiv:1608.08114}.
We do not have a consultation with specialists to comment the status and perspectives. 
\end{remark}

\subsubsection{Generalisations over $k$ and obstructions over $B$}\label{subsubsect:eneralisationsoverkandobstructionsoverB} 
The pioneer article \cite{Q} gave rise 
to a list of arguments 
that recover and generalise the result 
as the acyclicity of 
Cousin complexes
associated with some classes of cohomology theories over base fields,
which includes 
\begin{itemize}
\item
\cite{BO} by Bloch and Ogus 
concerned 
cohomology theories
that got their names,

\item
\cite{Voe-hty-inv} 
by Voevodsky on pretheories with transfers over perfect $k$
\footnote{ 
extended later to K-transfers 
\cite{KnizelSHIKor} 
and
framed transfers 
\cite{hty-inv}, 
and extended to arbitrary base fields in \cite{StrictA1invariancefields}
},
\item
\cite{Morel-connectivity,M3} by Morel 
for  
$\SHSone(k)$-representable coh. th.%
\footnote{
using Gabber's Presentation Lemma 
\cite{Gab,CTHK} for infinite $k$, 
proven later for finite fields in \cite{HK}
},
\item
\cite{Pan19} by Panin for 
coh. th. 
is sense of 
\cite[Def. 2.1]{PanafterSmirnovOrCoh}
for any $k$.
\end{itemize}
\par Summarising, 
\begin{itemize}[leftmargin=28pt]
\item[]
\begin{tabular}{||p{15cm}}
the 
$\A^1$-invariance 
and
Nisnevich excision
of 
a 
cohomology theory $A$ on $\Sm_k$
are
enough 
for 
the acyclicity of the associated Cousin complex 
over any base field $k$. 
\end{tabular}
\end{itemize}

\begin{remark}[Stable Connectivity over $B$]\label{sect:Stable-Connectivity-Counter-remark-over-B}

The acyclicity of the Cousin complexes
for all $\SH(B)$-representable theories
on the schemes $U=X^h_x$, 
$X\in \Sm_B$,
would 
imply 
the statement of 
Connectivity Theorem for $\SH(B)$ 
by the arguments of \cite{Morel-connectivity}
for any noetherian separated base scheme $B$ of finite Kull dimension.
The latter claim got name 
Connectivity Conjecture
and 
was
disproven in
\cite{AyoubcontrexempleA1connexite} by 
Ayoub. 
\end{remark}

\begin{remark}[Counter-example over $B$]\label{sect:Counter-remark-over-B}
The acyclcity 
fails for 
the cohomology theory represented by $\Sigma^\infty_{\PP^1}(B-z)$ for any 
closed point $z$ of positive codimension of any $B\in\Sch$, 
\cite{ColumnsCousinbiCompexSHfrB}.
So
\begin{itemize}[leftmargin=28pt]
\item[]
\begin{tabular}{||p{15cm}}
the 
$\A^1$-invariance and Nisnevich excision, and 
even $\SH(B)$-representability,
are not enough for the acyclcity of the Cousin complex over $B$, whenever $\dim B>0$,
\end{tabular}
\end{itemize}
and the proof for \eqref{eq:KthComplexUn}
should 
use other structures/poperties on/of the K-theory.
\end{remark}

\subsubsection{Results for 
``parts or modifications'' of K-theory}
\label{sect:ProofsforsomepartsormodificationsofK-theoryoverB}
T%
he acyclicity of 
\eqref{eq:KthComplexUn} 
is proven 
for
\begin{itemize}
\item
K-theory with 
finite 
coefficients
$K_*(-,\mathbb Z/n)$
\cite{zbMATH03955070} by Gillet 
out of 
the characteristic, and 
\cite{zbMATH01443583} by 
Geisser 
and 
Levine 
at the characteristic,
\item
Milnor K-theory 
$K^\mathrm{M}_*(-)$
%
\cite{zbMATH07659868}
by 
L{\"u}ders 
and 
Morrow, 
%
\end{itemize}


\begin{remark}[Surjectivity Property]
%
For 
any 
coh. th. $A$ 
defined by 
\begin{itemize}
\item[+] functors $(X,X-Z)\mapsto A^*(X,X-Z)=:A^*_Z(X)$ for $X\in \Sch_B$, and
\item[+] differentials $\partial\colon A^*(X-Z)\to A^{*+1}_{Z}(X)$, 
\end{itemize}
and
equipped additionally
with 
\begin{itemize}
\item[+]
a module structure 
over the Milnor-Witt K-theory compatible with 
the differentials, 
\item[+]
Gysin isomorphisms $A^*(x)\cong A_{x}^{*}(X)$, $x\in X$, 
compatible with the previous points,
\end{itemize}
the acycliclty of 
the Gersten complex for $A$ on 
any DVR spectrum $U\in \Sch_B$
can be deduced from 
the surjectivity 
of the restriction homomorphism 
\begin{equation}\label{eq:AXhxtwoheadarrowAx}A^*(U)\twoheadrightarrow A^*(x),\quad x\in U^{(1)},\end{equation}
which 
holds 
(0) %
for $U$ equipped with a retraction to $x$
for any $A$, 
and
(1)
for 
any 
local $U$
for
\begin{itemize}
\item 
$K^\mathrm{M}_*(-)$
by lifting of the generators,
and

\item 
$K_n(-)$, $n\leq 2$,
because of the isomorphism $K^\mathrm{M}_n(x)\to K_n(x)$,

\end{itemize}
and 
(2) for
Henselian local 
$U$
for
\begin{itemize}

\item 

$K_*(-,\mathbb Z/n)$
by 
the Rigidity Theorem 
\cite{SuslinKthlocalfields,GabberKthHenselianpairs} by Suslin and Gabber, 
and

\item 
$K_*(-)$
for 
such $U$
that
the 
residue field at $x$ 
being 
an algebraic extension of a finite field
because of the computation 
of K-theory of finite fields  
\cite{Quillen_K-thfinitefield} by Quillen, 
\end{itemize}
and 
it would follow for  
\begin{itemize}
\item 
$K_3(-)$
for any 
Henselian local 
$U$
from
the Rigidity Conjecture \cite[Conjecture 11.7]{zbMATH04168919}%
.
\footnote{
The
latter two points are mentioned to the author by A. Ananyevsky.
}
\end{itemize}
Nevertheless,
\begin{itemize}[leftmargin=30pt]
\item[]
\begin{tabular}{|p{15cm}}
it is
unknown 
do 
the surjectivity \eqref{eq:AXhxtwoheadarrowAx} 
hold
for 
higher groups
$K_*(-,\mathbb Q)$ 
for
Henselian DVR spectra $U$ of 
a characteristic 
$(0,p)$, $p>0$,
with 
an 
infinite residue field 
at $x$.
\end{tabular}
\end{itemize}

\end{remark}

\subsubsection{Reduction to one-dimensional schemes}\label{subsubsect:Reductiontoone-dimensionalschemes} 
The acyclicity of \eqref{eq:KthComplexUn} 
for 
essentially smooth one-dimensional 
schemes $U$ over $\mathbb Z$
implies 
the acyclicity of \eqref{eq:KthComplexUn} 
\begin{itemize}
\item 
for all schemes of the form $U=X_x$
\cite{zbMATH04122110} by 
Gillet 
and
Levine, 
and
\item 
for all regular unramified local rings 
\cite{arXiv:1710.00303} by Skalit.
\end{itemize}

\subsubsection{Reduction to the injectivity claim}\label{subsubsect:Reductiontotheinjectivityclaim} 
\cite{MR0862630} by 
Bloch 
provides a quasi-isomorphism of 
the complex \eqref{eq:KthComplexUn}
and the 3-term complex 
\begin{equation}\label{eq:KthtfComplexUn}K_n(U)\to K_n(U-U\times_B z)\to K_{n-1}(U\times_B z)\end{equation}
%
where $z\in B$ is the image of the closed point $x\in U$,
for any essentially smooth local $U$ scheme over a DVR.\footnote{Generalised 
by \cite{arXiv:1710.00303} by Skalit 
to geometrically regular $U$ over DVR
.}
\footnote{Generalised by \cite{StrictA1invariancefields,ColumnsCousinbiCompexSHfrB} to $\SH(B)$-representable coh. th. 
for any 
$B\in\Sch$, $\dim B<\infty$.
}
Consequently, 
\begin{itemize}
\item[-]
the complex \eqref{eq:KthComplexUn} is acyclic 
except the terms $K_n(U)$ and $\oplus_{x\in U^{(1)}}K_{n-1}(x)$%
,
\footnote{
The bounded acyclicity 
is proven also in \cite{arXiv:1710.00303} 
for geometrically regular Henselian $U$ over $B\in\Sch$.%
}
\footnote{
The 
acyclicity 
except the mentioned terms
is proven 
also 
in \cite{arXiv:2202.00896} by Panin
for semi-local schemes $U$ of 
smooth schemes
and
coh. th. in sense of \cite[Def. 2.1]{PanafterSmirnovOrCoh} over a DVR.}
and 
\item[-]
the acyclicity of \eqref{eq:KthComplexUn} 
reduces
to the injectivity 
$K_*(U)\hookrightarrow K_*(U-U\times_B z)$.
\end{itemize}
%
%

\subsection{Proof for essentially smooth local Henselian schemes} 
Constructing 
of 
Nisnevich squares and $\A^1$-homotopies 
over 
netherian separated schemes $B$ of finite Krull dimension 
allowed 
in \cite{SS,DHKY,zbMATH07612787}
proving of the triviality of 
the support extension morphisms
\begin{equation}\label{eq:extensionofsupport}
A_W(U)\to A(U)
\end{equation}
for any closed immersion $W\to U=X^h_x$ of positive relative codimension
for all 
pointed functors $A\colon \HpointedANis(B)\to \Set^\bullet$.
A 
novel 
instrument 
allows us to prove 
the claim for the closed fibre immersion 
\begin{equation}\label{eq:closedfibreimmersion}W=U\times_B z\to U=X^h_x
,\end{equation} 
where $z\in B$ is the image of $x\in X\in\Sm_B$, 
for motivic homotopies 
of cellular motivic $B$-spaces. 
\begin{nonum-theorem}[Theorems \ref{th:RLP:HsupK}, \ref{th:Injectivity}]

Let 
$B\in \Sch$ be regular of Krull dimension one.
Let $U$ be an essentially smooth local henselian scheme over 
$B$.
Let $\rho$ denote a regular function on $B$ that vanishing locus equals a closed point in $B$
as well as its inverse image on $U$.

For any non-negative integer $n$,
the morphism of the pointed sets
\begin{equation}\label{eq-intro:morphism:HsupKUZ->HsupKU}
\supextensionmorphismHAtopologyB
\colon 
\MapSet_{\HpointedAtopology(B)}(U_+/(U-Z(\rho))_+\wedge S^n,\motSpt)\to \MapSet_{\HpointedANis(B)}(U_+\wedge S^n,\motSpt)
\end{equation} 
is trivial
for any $\motSpc\in \HpointedAtopology(B)$ generated by 
$\PP^{\wedge l}_B$ for all $l\in \mathbb Z_{\geq 0}$ 
via colimits and extensions,
and
$\topology = \Zar$.
\end{nonum-theorem}
Though 
the result 
holds 
for the Nisnevich topology, 
which follows from 
\Cref{rem:BlrhocompintAoneABHBCBlrasNis-motivichomotopycategory:Injectivity},
for simplicity, in \Cref{sect:Presentation,sect:InvarianceLiftingPorpertyandInjectivity}, 
we present the proof 
for 
\[\topology = \Zar.\]
This is enough 
%
%
for
the triviality of \eqref{eq:extensionofsupport} 
for K-theory proven in \Cref{sect:Acyclicity},
because 
by \cite{zbMATH06773295}
presheaves of vector bundles on affine schemes are represented by Grassmanians in both $\Hpointed_{\A^1,\Zar}(B)$ and $\Hpointed_{\A^1,\Nis}(B)$. 
Then 
the acyclicity of 
the $3$-term complex \eqref{eq:KthtfComplexUn}
follows.
Since 
as mentioned above 
the complexes \eqref{eq:KthComplexUn} and \eqref{eq:KthtfComplexUn}
are quasi-isomorphic \cite{MR0862630}, and also \cite{arXiv:1710.00303,StrictA1invariancefields,ColumnsCousinbiCompexSHfrB},
\eqref{eq:KthComplexUn} is acyclic. 

The method for the triviality of \eqref{eq-intro:morphism:HsupKUZ->HsupKU} in \Cref{th:Injectivity} 
in our article is 
the lifting 
of the elements 
of \[\MapSet_{\SHSone_{\A^1,\topology}(B)}(U_+/(U-Z(\rho))_+\wedge S^n,\motSpt)\]
along 
the closed immersion of open pairs
\[(U,U-Z(\rho)) \xrightarrow{\rhomap} (\A^N_{U},\A^N_{U}-Z(t_0))^{h}_{\rhomap(U)},\]
where $\rhomap\colon U\to \A^1_U$ is the morphism defined by
$\rho$%
, 
Theorem \ref{th:RLP:HsupK},
which reduces 
the claim
on \eqref{eq-intro:morphism:HsupKUZ->HsupKU}
to 
the mentioned above one 
on 
\eqref{eq:extensionofsupport} for 
$\codim_{U/B}W>0$, cited as \Cref{th:Injectivitypositiverelativecodimansion}.
The proof of the latter lifting 
in \Cref{sect:InvarianceLiftingPorpertyandInjectivity}%
, \Cref{th:RLP:HsupK},
uses two ingredients:
%
%
\begin{itemize}
\item[(1)] 
a new type of motivic equivalences called $\phi$-motivic equivalences 
defined as morphisms in a certain $\infty$-category denoted $\catOB$ 
in 
\Cref{section:ABHBCmotiviclocalisation}
based on \cite{ABHBCinvariance,ABHBClocalisation},  
\item[(2)]
an
interpretation of the morphisms 
$U/(U-Z(\rho))\wedge S^l\to\PP^{\wedge n}_B$ 
in $\HpointedAtopology(B)$
in 
terms of 
the category of presheaves $\Spc(\catOB)=\Func(\catOB^\mathrm{op},\Spc)$ 
in 
\Cref{sect:Presentation} based on \cite{cciblowupsrho}. 
    \end{itemize}
\par 
\subsubsection{$\phi$-motivic localisation and 
homotopy category}
The localisation of
the $\infty$-category $\Spc(\catOB)$ 
with respect to
the $\phi$-motivic equivalences 
gives rise further to 
the 
$\phi$-motivic homotopy category over $B$, 
\Cref{def:BlrhocompintAoneABHBCClrasequivtopology-motivichomooptycategory}.
%
The 
triviality like for 
\eqref{eq-intro:morphism:HsupKUZ->HsupKU}
holds for any 
pointed presheaf 
defined on the $\phi$-motivic homotopy category over $B$, \Cref{rem:BlrhocompintAoneABHBCBlrasNis-motivichomotopycategory:Injectivity},
and
the acyclicity of the Cousin complexe 
hold for any $\phi$-motivic prehseaf of spectra.
\begin{remark}\label{rem:ConnectivityandStrictHomooptyInvariance}
To continue the discussion from \Cref{sect:Counter-remark-over-B},
we note that
following further
the arguments of \cite{Morel-connectivity} 
we can deduce Stable Connectivity and Strict Homotopy Invariance Theorems
for the $S^1-$ and $\PP^1$-stable $\phi$-motivic categories
from the acyclcity of the Cousin complexes.
\end{remark}

The role played by the category $\Sm_B$ in the Morel-Voevodsky construction 
is played by 
the 
$\infty$-category $\catOB$, 
which definition   
in \Cref{sect:catAOB,sect:MorphismsTypes,sect:catOB},
starts with the category $\precatAOuB$ of 
data
\begin{equation}\label{eq:(schX(rhoSet,gr,supp))}
(\schX,(\rhoSet,(\gr,\supp,\presetsSchrhoBl)))
,\end{equation}
where 
$\schX\in\Sch_B$, $\rhoSet$ is some set of sections of line bundles,
$\gr,\supp$ are $\rhoSet$-monomials,
and $\beta$ is some set of sets of $\rhoSet$-monomials. 
The interpretation is as follows:
\begin{itemize}
\item $\schX$ 
is the underlying scheme of 
a compacitification,
\item $\gr$ relates to the ``infinitely remote subscheme'' called 
boundary of the compactification,
\item $\supp$ relates to the support of cohomologies,
\item $\presetsSchrhoBlcup$ define centres of some 
blow-ups on $\schX$,
\item $\rhoSet$ are parameters on $\schX$ that control $\gr$, $\supp$ and $\presetsSchrhoBlcup$.
\end{itemize}
%
%
\par
Then we define 
$\catOB$ 
as
certain full subcategory of 
the $\infty$-category $\Spc^{\precatAOB}$ 
of copresheaves on 
the category of arrows $\precatAOB = (\precatAOuB)^{\Delta^1}$,
\Cref{def:catOB}. 
The conditions on the objects \begin{equation}\label{eq:oxincatOB}\aroX = (\oX\to\oXprime)\in\catOB
, \quad \oX\to\oXprime\in \Spc^{\precatAOuB}
,
\end{equation}
in sense of functor
$\catOB\hookrightarrow 
\Spc^{\precatAOB}\simeq (\Spc^{\precatAOuB})^{\Delta^1}$,
imply 
the diagram 
in the $\infty$-category of copresheaves on $\Sch_B$
\[
\prod_{\alpha}\PP^{N_\alpha}_B\simeq 
\mathsf{P}
\xleftarrow{b}\limBlPobject
\xleftarrow{e}
\limEBlPobject
 \simeq\schX
\]
where
\begin{itemize}
\item
$b$ is a limit of blow-ups that centres are intersections of coordinate divisors and exceptional divisors,
\item
$e$ is a limit of affine \'etale morphisms,
\end{itemize}
%
\par
The $\phi$-motivic homotopy category
is 
the localisation $\Spc_{\tau,\arrowequiv,\Blrasequiv,\ABHBC,\compintA,\Blrho}(\catOB)$ of
the $\infty$-category of presheaves $\Spc(\catOB)$ 
%
%
%
with respect to
\begin{itemize}
\item $\tau$-equivalences, where $\tau$ is the topology on $\precatAOuB$ induced by a topology $\tau$ on $\Sch_B$,
\item
$\supp$-equivalences equivalence $\oX$ 
with 
$\operatorname{cofib}((\schX,(\rho,(\gr\cdot\supp,1,\presetsSchrhoBl)))\to (\schX,(\rho,(\gr,1,\presetsSchrhoBl))))$
\item
$\arrowequiv $-equivalences defined by the counit of the 
adjunction 
$\precatAOB\simeq (\precatAOuB)^{\Delta^1} \leftrightarrows \precatAOuB$,
\item
$\Blrasequiv$-equivalences, which allow modifying $\rhoSetprime$ on $\schXprime$ in \eqref{eq:oxincatOB},
\item
$\Blrho$-equivalences, which are 
limit morphisms of
certain 
blow-ups 
controlled by $\rhoSet$%
, \Cref{sect:Blbeta-localisation}, 
\item
$\compintAone$-equivalences, which are like $\square$-equivalences in \cite{MotiveswithmodulusIII,binda2024logarithmicmotivichomotopytheory},
%
%
\item
$\ABHBC$-equivalences%
,
which are the main ingredient for the proof of 
the lifting property 
of 
elements in motivic homotopy groups
w.r.t. 
closed immersions of essentially smooth local henselian $B$-schemes.
\end{itemize}

\begin{remark}[Illuminating 
$(\arrowequiv,\Blrasequiv,\supp)$-equivalences
and the data $\rhoSet$%
]\label{rem-intro:ComparisonTheoremilluminatingarrowequivequivalencesandrhoSet}
The $\phi$-motivic homotopy category can be constructed 
equivalently 
skipping
the data $\rho$ and $\supp$ 
and the arrows category $\precatAOB$ 
in the definitions of $\precatAOuB$ and $\catOB$,
\Cref{eq:ComparisonTheorem}.
Nevertheless, all of this is useful for 
the proof of the triviality of support extension maps and acyclcity of Cousin complexes
discussed in the article.
\end{remark}

\begin{remark}[
Continuous presheaves without $\ABHBC$-invariance condition]
There is 
a relation of 
the 
subcategory of continuous objects
in the 
category of $(\tau,\arrowequiv,\Blrasequiv,\compintA,\Blrho)$-invariant objects,
i.e. without $\ABHBC$-invariance condition,
to 
the $\square$-homotopy motivic categories form 
\cite{MotiveswithmodulusIII,binda2024logarithmicmotivichomotopytheory}. 
A 
difference of 
the localisation procedures,
nevertheless, 
%
is the type of blow-ups.
Resolution of singularities over $B$
could help to illuminate this difference. 

\end{remark}
\subsubsection{Interpretation of Morel-Voevodsky $\A^1$-homotopies in the $\infty$-category $\catOB$} 

The Part (2) 
interprets 
elements of 
\begin{equation}\label{eq:MapZarAUqUmZPnB}\Map_{\Hpointed_{\A^1,\Zar}(B)}(U_+/(U-Z(\rho))_+\wedge K,\PP^{\wedge n}_B)\end{equation}
as morphisms 
in $\Spc_{\Blrho,\Blrasequiv}(\catOB)$
from the object in $\catOB$ defined by \[(U,(\rho,1,\rho))\in\precatAOuB\] 
to 
the motivic sphere $\PP^l_B$ in $\precatAOuB$.
The construction starts with 
coding of the element in \eqref{eq:MapZarAUqUmZPnB} 
by
sets of sections of line bundles on 
the 
compactified geometric realisation $|K|_{\compintA}\times U$ of $K$ over $U$.
%
The latter sections are replaced by sections on certain blow-up of $|K|_{\compintA}\times U$ using \Cref{cor:surj:sinLonBlcXcapZsiseqmpty->sinLonXcapZssubsetcupZs} from  
\Cref{section:subsectionPresentation}
based on the 
study of blow-ups of schemes 
with 
relation to the data $\rhoSet$, \Cref{def:precatAOuobjects} 
applied 
to blow-ups in the category $\catOB$%
.



\subsection{Awknowlagement}
The author
is thankful to I.~Panin,
who involved the author to the theme.
The author is thankful to A. Kuznetsov for a set of series of consultations during the research.
The author thanks also M. Temkin and A. Ananyevsky. 
\subsection{Notation}

\begin{enumerate}

\item
$\codim_X(Z)$ is the codimension of $Z$ in $X$.
$\codim_X(\emptyset)=\infty$. 

\item
$\codim_{U/B}W$ is the relative codimension of $W$ in $U$ over $B$.

\item
$\MapSet_C(-,-)$, $\MapSpc_C(-,-)$ and $\MapSpt_C(-,-)$
denote respectively the mapping set or the mapping space, or the mapping spectrum in 
a category or an $\infty$-category $C$, or a stable $\infty$-category $C$.

\item
$S^1_{\A^1} \simeq \Delta^1_{\A^1}/\partial \Delta^1_{\A^1}$,
$S^l_{\A^1} \simeq (S^1_{\A^1})^{\wedge l}$

\item
Let $\ShiftsColimits_{C}(-)$ denote the subcategory of $C$ generated via shifts and colimits.

\item
$\Subsets(S)$ is the set of subsets of a set $S$.

\item\label{notation:I(Z)_AND_Z(I)}
Let $X$ be a scheme.
\par $X^{(c)}$
is the set of points on codimension $c$.
\par
$X_x$ and $X^h_x$ are the local and Henselian local schemes of a scheme $X$ at $x\in X$.
\par
Given set of sections of line bundles $I$ on $X$,
$Z(I)=Z_X(I)$ is the vanishing locus.
\par
Given closed subscheme $Z$ in $X$,
$I(Z)$ is the sheaf of vanishing ideals.
\par
Given divisor $D$ on $X$,
$I(D)=I_X(D)$ is the sheaf of ideals defined by $D$.
Let \[Z_X(D):=Z_X(I_X(D)).\]



\item
\label{notation:productcdotinMon}
$A\cdot B=\{a\cdot b\in \Monoid| a\in A, b\in B\}\in \Subsets(\Monoid)$ 
for some $A,B\in \Subsets(\Monoid)$.
%

\item\label{notation:setofsubsets}
$2^S$ is a set of subsets of a given set $S$.
\item\label{notation:MonS}
$\Mon(S)=S^\Mon$ is the monoid of monomials of a given set $S$.\\
\item\label{notation:MonexpS}
We use the identification of $\Mon(2^S)=(2^S)^\Mon$ and the submonoid is the monoid of subsets of $S^\Mon$ generated by subsets of $S$, Notation \ref{notation:productcdotinMon}.

\item\label{notation:setsrhoBlOperations}
A $\setsrhoBlOperations$-closed set 
is a set 
$\beta$ such that
for each $\beta_1,\beta_2\in \beta$, $\beta_1\cup \beta_2\in \beta$.
$\beta^{\setsrhoBlOperations}$ 
is the minimal $\setsrhoBlOperations$-closed subset containing given subset $\beta$.
%
%

\item\label{notation:Spcpointed}
$\Spcpointed$ denotes the $\infty$-category of pointed spaces.
\item\label{notation:SuC_AND_S(C)}
$S^C=\Func(C,S)$,
$S(C)=\Func(C^\op,S)$.
\item
$\Sch$ and $\Sch_B$ are the categories of Noetherian separated schemes and $B$-schemes.

\item 
$\Spc$, $\Spt$ are $\infty$-categories of spaces and spetra.
\item 
$\Spc_{\A^1,\tau}(\Sch_B)$ and $\Spt_{\A^1,\tau}(\Sch_B)$
are the subcategories of $\A^1$-invariant Nisnevich excisive presheaves of spaces and spectra on $\Sch_B$, and similarly for $\Sm_B$.
$\Spt_{\A^1,\tau}(B)=\Spt_{\A^1,\tau}(\Sch_B)$.


\item
\label{def:continuous}
Given $\infty$-category $C$,
a presheaf $F\in\Spc(\Spc^C)$ is continuous, if 
for any $X \simeq \varprojlim X_\alpha$ in $C$,
the canonical morphism 
$F(X)\leftarrow \varinjlim_{\alpha} F(X_\alpha)$
is an equivalence.

\item
\label{notation:associatedcontinuous}
For any presheaf $F\in \Spc(C)$ on an $\infty$-category $C$,
the same symbol $F$ denotes
the associated continuous presheaf on the $\infty$-category of copresheaves $C^{\Spc}$, Notation {def:continuous}.

\item \label{def:categoryDeltaleqone}
$\Delta^\bullet$ denotes the category of finite linearly ordered sets with objects 
$\{0,\dots n\}=[n]=\Delta^n$.\\
$\Delta^{\bullet\leq l}$ is the subcategory spanned by $\Delta^{n}$, $n\leq l$.\\
$\Deltatimesbullet=(\Deltableqone)^{\mathbb{Z}_{\geq 0}}$.

\item \label{def:cubicalset}
A \cubicalset $S$
is a functor
\[S\colon (\Deltatimesb)^\mathrm{op}\to \Set.\]

\item\label{def:Mon}
Let $\rho$ be a subset in a monoid $N$,
$\Mon(\rho)$ denotes the submonoid generated by $\rho$.

\item\label{notation:Sch_AND_SmB}
$\Sch$ is the category of noetherian separates schemes.\\
$\Sch_B$ is the category of $B$-schemes of such type.\\
$\Sm_B$ is the category of smooth $B$-schemes.\\


\item
\label{notation:LLP_AND_RLP}
Let 
$c\colon X\to Y$, 
$f\colon F\to G$
be morphisms of an $\infty$-category.
The morphism $c$ satisfies LLP w.r.t. the morphism 
$f$
or
$f$ satisfies RLP w.r.t. $c$
if
for any morphisms $X\to F$ and from $Y\to G$ such that the outer square in the diagram is commutative
\[\xymatrix{X\ar[d]\ar[r] & F\ar[d]\\Y\ar[r]\ar@{-->}[ru]&G}\]
the space of dashed arrows that fit into the commutative diagram
is non empty.

\item
\label{notation:uLLP_AND_uRLP}


Let 
$c\colon X\to Y$, 
$f\colon F\to G$
be morphisms of an $\infty$-category.
The morphism $c$ satisfies uLLP w.r.t. the morphism 
$f$
or
$f$ satisfies uRLP w.r.t. $c$
if
for any morphisms $X\to F$ and from $Y\to G$ such that the outer square in the diagram is commutative
\[\xymatrix{X\ar[d]\ar[r] & F\ar[d]\\Y\ar[r]\ar@{-->}[ru]&G}\]
the space of dashed arrows that fit into the commutative diagram
is contractible,
i.e.
the morphism
\[
\Map(Y,F)
\rightarrow
\Map(X,F)\times_{\Map(X,G)}\Map(Y,G).\]

%
%

\item\label{notation:hens}
Given a morphism of schemes $X\to X^\prime$,
and closed subscheme $Z^\prime$ in $X^\prime$,
$X^h_{Z^\prime}:=X^h_{Z^\prime\times_{X^\prime} X}$. 

\item\label{notation:setminuszerocalS(X)}
$\setminuszero{\calS}(X)=\calS(X)\setminus \{0\}=\calS(X)-0$.

\item\label{notation:schxus(setminuszerocalS)}
Given $f\colon X\to X^\prime$ in $\Sch$, 
$f^*\colon \ssetminuszero{\calS}(X^\prime)\dashrightarrow \ssetminuszero{\calS}(X)$
is a partly defined map
induced by $f^*\colon \calS(X^\prime)\to \calS(X)$ and immersion $\setminuszero{\calS}\to\calS$.


\end{enumerate}

\section{$\ABHBC$-motivic localisation}
\label{section:ABHBCmotiviclocalisation}

\subsection{Schemes category preparations}
\subsubsection{Line bundles and sections} 
\begin{definition}\label{def:LIStoLtoT}
Define the morphism of presheaves
on $\Sch$
\begin{equation}\label{eq:LIStoL}
\calS\to \calL; \quad
(s,L)\mapsto L,
\end{equation}
where
$\calL$ denotes the presheaf of 
line bundles on $\Sch$,
and 
$\calS$ is the presheaf of line bundles sections,
i.e.
\[\calS(X)= \{(s\in L(X),L\in \calL(X))\}.\]
\end{definition}


\begin{remark}\label{rem:MonoidcalS}
The presheaf of sets $\calS(X)$ has the monoid structure 
with respect to the tensor products of line bundles and products of sections.
\end{remark}
\begin{remark}\label{rem:MonoidcalS:invetibleelements}
Each submomoid of $\calS(X)$ contains all invertible sections, i.e. elements $(s,L)$ such that $Z_X(s)=\emptyset$, because the set $\calS(X)$ is the set of isomorphicm classes of the pairs $(s,L)$.
\end{remark}

\begin{definition}\label{def:irreduciblesectionoflinebundle}

Given $X\in \Sch$, an element $s\in \calS(X)$ is \emph{irreducible} if the scheme $Z_X(s)$ is irreducible.
\end{definition}

\subsubsection{Subpresheaves of $\calS^{\times l}$}

\begin{definition}\label{def:calSpowerSoverLcalL}
We use notation $\calS^{\times l}_{\calL}$ for powers of $\calS$ over $\calL$ w.r.t. the morphism \eqref{eq:LIStoL},
and
\[\calS^{\times\powerSoverL}_{\calL}:=\prod_{l\in \powerL}\calS^{\times\powerS\times_{\powerL}\{l\}}_\calL\]
for a given morphism of sets
$\powerSoverL\colon \powerS\to\powerL$.

\end{definition}

\begin{definition}\label{def:sinLcapZsequalsZdelta}
The category $\sinLcapZsequalsZdelta$ is a full subcategory in $\Set(\Sch)$
that objects are presheaves $F$
of the form 
\begin{equation}\label{eq:sLonXKZsequalsZdelta} 
F(X) = 
\{
s\in \calS^{\times\powerSoverL}_\calL
| 
\cap_{\indexSubsetssrho\in \SetofindexesSubsetsrho_\indexSubsets}Z( s_{\indexSubsetssrho} )
\subset 
Z(\delta_{\indexSubsets}),
\indexSubsets\in \SetofindexesSubsets
\}
\end{equation}
for 
a morphism of finite sets $\powerSoverL\colon \powerS\to\powerL$
and
a set $\SetofindexSubSets$
with
subsets $\SetofindexesSubsetsrho_\indexSubsets\subset \powerS$
and 
sections
$\delta_{\indexSubsets}\in \calS(X)$
for each $\indexSubsets\in \SetofindexSubSets$.
\end{definition}

\begin{definition}\label{def:sinLcapZsequalsZempty}
The category $\sinLcapZsequalsZeqmpty$ is a full subcategory in $\Set(\Sch)$
that objects are presheaves $F$
of the form 
\begin{equation}\label{eq:sLonXKZsequalsZempty}
F(X) = 
\{
s\in \calS^{\times\powerSoverL}_\calL
| 
\cap_{\indexSubsetssrho\in \SetofindexesSubsetsrho_\indexSubsets}Z( s_{\indexSubsetssrho} )
\subset 
\emptyset,
\indexSubsets\in \SetofindexesSubsets
\}
\end{equation}
for 
a morphism of finite sets $\powerSoverL\colon \powerS\to\powerL$
and
a set $\SetofindexSubSets$
with
subsets $\SetofindexesSubsetsrho_\indexSubsets\subset \powerS$
for each $\indexSubsets\in \SetofindexSubSets$.
\end{definition}

\subsubsection{Retract closed immersions}

\begin{definition}\label{def:copresheavesSch:retractclosedimmersion}
An morphism $\climmersion\colon X\to\ovX$ in $\catAOuB$ is a 
retract-closed immersion
if
it is a closed immersion and
there is a morphism
$\ovX\to X$
such that the composite 
\[X\xrightarrow{\climmersion}\ovX\to X\]
equals $\mathrm{id}_{X}$.
\end{definition}

\subsubsection{Equational morphisms}

\begin{definition}\label{def:Sch:byequations}
A morphism $X\to Y$ in $\Sch$ is \emph{\byequations }
if
it is isomorphic to the composite
\[X\simeq Z_{\A^N_Y}(\eX) - Z_{\A^N_Y}(\ecompltoX)\to \A^N_Y\to Y\] 
for some sets of regular functions $\eX$ and $\ecompltoX$ on $\A^N_Y$.
\end{definition}

\subsection{$\infty$-category $\catAOB$} 
\label{sect:catAOB}

\subsubsection{Category $\precatAOuB$} 

\begin{definition}\label{def:precatAOuobjects}
An $\precatAOu$-object over $B$ 
is the 
data
\begin{equation}\label{eq:oX}\oX=(\schX,(\rhoSet,(\gr,\supp,\presetsSchrhoBl))),\end{equation} 
where
\begin{itemize}
\item[(o1)] 
$\schX\in \Sch_B$, 
\item[(o2)] 
$\rhoSet\in \Subsets(\setminuszero{\calS(X)})$, Notation \ref{notation:schxus(setminuszerocalS)}, 
such that for each finite subset $\rhoSubset\subset\rhoSet$,\begin{itemize}
\item
$\codim_\schX(Z(\rhoSubset))\geq\#\rhoSubset$,
\item
connected components of the scheme $Z_X(\rhoSubset)$ are irreducible and reduced,
and moreover, \item the scheme $Z_X(\rhoSubset)$ is irreducible if $\#\rhoSubset=1$, i.e. 
for each $\rhoElement\in \rhoSet$, the scheme $Z_X(\rhoElement)$ is irreducible,
\end{itemize}
\item[(o3)] $\gr,\supp\in \Monrho$, where $\Monrho=\Mon(\rhoSet)$ is the submonoid of $\calS(X)$ generated by $\rhoSet$, \Cref{rem:MonoidcalS},

\item[(o4)] 
$\presetsSchrhoBlcup$ is a $\cup$-closed set of finite subsets of 
$\Monrho$, Notation \ref{notation:setsrhoBlOperations}, \ref{notation:MonexpS}. 
\end{itemize}
\end{definition}

\begin{notation}\label{notation:presetsSchrhoBlcupANDsetsSchrhoBl}
Denote by $\setsSchrhoBl$ the submonoid in 
the monoid of subsets of the monoid of monomials $\Monrho$
generated by $\presetsSchrhoBlcup$
. 
\end{notation}

\begin{definition}\label{def:precatAOumorphisms}
An $\precatAOu$-morphism over $B$
\begin{equation}\label{eq:oXtooXprime}
\ox\colon
\oX=(\schX,(\rhoSet,(\gr,\supp,\presetsSchrhoBl)))\to 
\oXprime=(\schXprime,(\rhoSetprime,(\grprime,\suppprime,\presetsSchrhoBlprime)))
\end{equation}
is defined by 
\begin{itemize}
\item[(m1)]
a morphism of $B$-schemes
\begin{equation}\label{eq:schXtoschXprime}\schx\colon \schX\rightarrow\schXprime\end{equation}
such that
\item[(m2)] 
$\Monrho \supset \schx^*(\Monrhoprime)$, 
\item[(m3)] 
$\gr\leq\schx^*(\grprime)$,
$\supp\leq \schx^*(\suppprime)$,
\item[(m4)] 
$\setsSchrhoBl\supset \schx^*(\setsSchrhoBlprime)$,
\end{itemize}
\end{definition}
\begin{remark}\label{rem:morphism:rhoSet}
For any morphism \eqref{eq:oXtooXprime} in $\precatAOuB$,
the set of subschemes $Z_X(\rhoSet)$ for all $\rhoSet\in \calS(\schX)$ 
contains 
the set of irreducible components of the subschemes 
$Z(\schx^*(\rhoSetprime))$ for all $\rhoSetprime\in \calS(\schXprime)$.
\end{remark}
\begin{definition}\label{def:precatAOu}
$\precatAOuB$
is the category 
that objects and morphisms are 
$A$-objects and $A$-morphisms over $B$, 
\Cref{def:precatAOuobjects,def:precatAOumorphisms},
and the composition of morphisms of the form $\ox$ \eqref{eq:oXtooXprime} 
is induced by the composition of the respective morphisms $\schx$ \eqref{eq:schXtoschXprime} 
in $\Sch$.
\end{definition}
%

%

\subsubsection{Category $\precatAOB$}
\begin{definition}\label{def:precatAO}
$\precatAOB$ is the category 
of arrows of $\precatAOuB$, i.e. the morphisms of the category $\precatAOuB$, i.e.
\[\precatAOB:=(\precatAOuB)^{\Delta^1}.\]
\end{definition}
\begin{notation}\label{notation:morphisms_precatAOB}
We use notation \begin{equation}\label{eq:vecoXcolonoXtooXprime}
\aroX\colon
\oX=(\schX,(\rhoSet,(\gr,\supp,\presetsSchrhoBl)))\to 
\oXprime=(\schXprime,(\rhoSetprime,(\grprime,\suppprime,\presetsSchrhoBlprime)))
\end{equation}
for objects of $\precatAOB$.
\end{notation}

\begin{definition}\label{def:morphisms:precatAOuB->precatAOB}
Let $\classa$ be a class of morphisms in $\precatAOuB$ closed with respect to base changes along arbitrary morphisms in $\precatAOuB$.
%
We define the class of morphisms $\classa$ in $\precatAOB$ as the class 
generated by 
the classes $\rasconstmorphismsclassa$ and $\rasmorphismsclassa$ via the composition.
Here
$\rasconstmorphismsclassa$ and $\rasmorphismsclassa$
denote
the classes of morphisms in $\precatAOB$
defined respectively by the commutative squares in $\precatAOuB$ of the form
\begin{equation}\label{eq:morphisms:precatAOuB->precatAOB}
\xymatrix{
\constroveroX\ar[r]^{\classa}\ar[d]\ar@{}[rd]|{\rasconstmorphismsclassa}& \oX\ar[d]\\ 
\constroveroXprime\ar[r]^{\simeq}_{\classa}& \oXprime
,}
\quad
\xymatrix{
\constroveroX\ar[r]^{\classa}\ar[d]\ar@{}[rd]|{\rasmorphismsclassa}& \oX\ar[d]\\
\constroXprime\ar[r]_{\classa}& \oXprime
,}
\end{equation}
where 
the vertical arrows define the objects in $\precatAOB$ and 
the horizontal ones define the morphisms,
\begin{itemize}
\item 
the left bottom arrow is an isomorphism, 
and 
the right square is a pullback,
and 
\item 
all the horizontal arrows morphisms belong to $\classa$.
\end{itemize} 



\end{definition}
%
%
%

\subsubsection{Monoidal structure}
%

\begin{definition}\label{def:precatAOuB:times}
Given $\oX_1,\oX_2\in \precatAOuB$,
\begin{multline*}\label{eq:precatOuB:times}
\oX_1\times\oX_2=
(\schX_1,\rhoSet_1,\gr_1,\supp_1,\presetsSchrhoBl_1)
\times
(\schX_2,\rhoSet_2,\gr_2,\supp_2,\presetsSchrhoBl_2)
:=\\
(\schX_1\times\schX_2,
\rhoSet_1\cup \rhoSet_2,
\gr_1+\gr_2,
\supp_1+\supp_2,
\presetsSchrhoBl_1\cup \presetsSchrhoBl_2
)
,\end{multline*}
where
$\rhoSet_\alpha$,
$\gr_\alpha$,
$\supp_\alpha$,
$\presetsSchrhoBl_\alpha$,
$\alpha=1,2$,
at the right side 
stand for the inverse images along the projections
$\schX_1\times\schX_2\to \schX_\alpha$.
\end{definition}


\subsubsection{Functors to schemes}

\begin{definition}\label{def:FunctorSch}
$\FschX$ is the functor 
\begin{equation}\label{eq:funct:precatOBtoSch:oX->schX}
\begin{array}{lclcl}
\FschX&\colon&\precatOB & \to & \Sch_B\\ 
&&(\ox\colon \oX\to\oXprime)&\mapsto& \schX,
\end{array}
\end{equation}
where $\oX$, $\schX$, and $\gr$ are from \eqref{eq:oXtooXprime}. 
\end{definition}

\begin{notation}
For any $F\in \Spc(\Sch_B)$,
the same symbol denotes
the inverse image of $F$ along the functor \eqref{eq:funct:precatOBtoSch:oX->schX}.
\end{notation}

\begin{definition}\label{def:funcXcomplinfX}
$\funcXcomplinfX$ is
the functor \begin{equation}\label{eq:funct:precatOBtoSch:X-Z(gr)}
\begin{array}{lclcl}
\funcXcomplinfX&\colon& \precatOB & \to & \Sch_B\\
&&(\ox\colon \oX\to\oXprime)&\mapsto& \schX-Z(\gr),
\end{array}
\end{equation}
where $\oX$, $\schX$, and $\gr$ are from \eqref{eq:oXtooXprime}. 
\end{definition}

%
\begin{definition}\label{def:funcXcomplinfXquotientXcomplinfXunionSup}
$\funcXcomplinfXquotientXcomplinfXunionSup$ is
the functor 
\begin{equation}\label{eq:funct:precatOBtoSch:(X-Z(gr))/(X-Z(gr)-Z(sup))}
\begin{array}{lclclcl}
\funcXcomplinfXquotientXcomplinfXunionSup&\colon& 
\precatOB & \to & \Spc(\Sch_B)\\
&&
(\ox\colon \oX\to\oXprime)&\mapsto& (\schX-Z(\gr))/(\schX-Z(\gr)-Z(\supp)),
\end{array}
\end{equation}
where $\oX$, $\schX$, and $\gr$ are from \eqref{eq:oXtooXprime}.
\end{definition}
The natural morphism in $\Spc(\Sch_B)$
\[(\schX-Z(\gr))\to (\schX-Z(\gr))/(\schX-Z(\gr)-Z(\supp))\]
defines the morphism of functors
\begin{equation}\label{eq:morphism:funcXcomplinfX->funcXcomplinfXquotientXcomplinfXunionSup}\funcXcomplinfX\to \funcXcomplinfXquotientXcomplinfXunionSup.\end{equation}

\subsubsection{$\infty$-category $\catAOB$}
%


\begin{definition}\label{def:catAO}
The $\infty$-category $\catAOB$ is the $\infty$-category ${(\precatAOB)}^{\Spc}$ of 
co-presheaves on the category of $\precatAOB$.
\end{definition}
\begin{notation}\label{notation:inverseimagealongthefucntoroncopresheaves_precatAOB->SchB}

For any $F\in \Spc(\Spc^{\Sch_B})$,
the same symbol denotes
both 
the inverse image along the functor \[(\FschX)^{\Spc}\colon(\precatAOuB)^{\Spc}\to (\Sch_B)^{\Spc}\] induced by \eqref{eq:funct:precatOBtoSch:oX->schX}
and
the inverse image along the composite with the functor 
\[(\FschX)^{\Spc}\colon(\precatAOB)^{\Spc}\to (\precatAOuB)^{\Spc}\] induced by 
the functor 
\begin{equation}\label{eq:functor:precatAOB->precatAOuB:(oXtooXprime)mapsto(oX)}\begin{array}{lcl}
\precatAOB&\to&\precatAOuB;\\
(\ox\colon \oX\to\oX^\prime)&\mapsto& \oX
.\end{array}
\end{equation}

\end{notation}%
%
%

\begin{remark}\label{rem:catAOB}
In view of the Grothendieck correspondence objects in $\catAOB$ are defined by diagrams in the category $\precatAOB$, which are diagrams in $\precatAOuB$ of the form 
\begin{equation}\label{eq:DtimesDeltaonetoprecatAOuB}D\times\Delta^1\to \precatAOuB,\end{equation}
where $D$ is an $\infty$-category indexing the diagram in $\precatAOB$.
\end{remark}
\begin{definition}\label{def:morphisms:precatAOB->catAOB}
Given class of morphisms $\classa$ in $\precatAOB$, the same symbol $\classa$ denotes the class of morphisms in $\catAOB$ that are defined by term-wise $\classa$-morphisms of diagrams in $\precatAOB$ according to \Cref{rem:catAOB}.
\end{definition}


%
%



\subsection{Types of morphisms in $\catAOB$.}
\label{sect:MorphismsTypes}
\subsubsection{Types of morphisms with respect to $\rhoSet$.} 

\begin{definition}\label{def:rhoirr_rhoimmersion}
A morphism 
$\ox\colon \oX\to \oXprime$
given by \eqref{eq:oXtooXprime}
is called 
\begin{itemize}
\item \emph{\rhoirr-morphism} if $\schx\colon \schX\to \schXprime$ induces a well defined map 
\[\rhoSetprime\to\rhoSet,\]

\item \emph{\preimagegrsup-morphism} if 
\[\gr = \schx^*(\grprime),\quad
\supp = \schx^*(\suppprime),\]

\item \emph{\preimagerho-morphism} if 
for any $\rhoelement\in \rhoSet$ there is $\rhoelementprime\in \rhoSetprime$
such that $Z_{\schX}(\rhoelement)\subset Z_{\schX}(\schx^*(\rhoelementprime))$,
\\ and 
$\gr = \schx^*(\grprime)$,
$\supp = \schx^*(\suppprime)$,
$\presetsSchrhoBlcup = \schx^*(\presetsSchrhoBlcupprime)$. 
\item \emph{\rhoimmersion-morphism} if 
there is a map 
$\rliftrho\colon \rhoSet\to \rhoSetprime
$
such that
\begin{gather*}
\schx^*(\rliftrho(\rhoelement))=\rhoelement,\; \forall \rhoelement\in\rhoSet,\\
\grprime = \rliftrho(\gr),\quad
\suppprime = \rliftrho(\supp),\\
\setsSchrhoBlprime = \{\rliftrho(\elementinsetsSchrhoBl)| \elementinsetsSchrhoBl\in \setsSchrhoBl\},
\end{gather*}
\item \emph{\rhotriv-morphism} if 
it is a 
\rhoimmersion-morphism 
such that
$\schx\colon \schX\to \schXprime$ induces a bijection 
\[
\rhoSetprime\xrightarrow{\simeq}\rhoSet
.
\]
\end{itemize}
\end{definition}

%

%
%
\begin{lemma}\label{lm:morphism:preimagerasrhogrsub}
For any $\oX\in \precatAOuB$ given by \eqref{eq:oX}
and morphism $\schx\colon \schX^\prime\to \schX$ in $\Sch_B$
,
the set of 
\preimagerasrho-morphisms $\ox\colon \oX^\prime\to \oX$ in $\precatAOuB$ given by \eqref{eq:oXtooXprime}
that image along the functor \eqref{eq:funct:precatOBtoSch:X-Z(gr)} equals 
$\schx$
is either empty or consists of unique element
.\end{lemma}
\begin{proof}
The inverse image along the morphism $\schx$ uniquely defines the data $(\rhoSet,\gr,supp,\presetsSchrhoBl)$.
If the latter data satisfies conditions form \Cref{def:precatAOuobjects} it defines the unique morphism $\ox\colon \oX^\prime\to \oX$, otherwise the set is empty.
\end{proof}

\begin{remark}\label{rem:rhoimmersion}
Any \rhoimmersion-morphism is a \rhopreimage-morphism.
\end{remark}

\begin{remark}\label{rem:rhotriv}
Any \rhotriv-morphism is a \rhoirr-morphism.
\end{remark}



\subsubsection{Closed immersions}

\begin{definition}\label{def:precatAOuB:rhoci-closedimmersion}
A morphism $\climmersion\colon \oX\to \ovoX$ in 
$\precatAOuB$ 
is a \emph{\rhoci-closed immersion}
if 
\begin{itemize}
\item
the morphism $\FschX(\climmersion)$ 
is a 
closed immersion in $\Sch$,
\item
there is a 
subset $\eqXinovX$ of 
$\Monrho(\ovoX)$
such that
\[\schX=Z_{\schovX}(\eqXinovX)
.\]
\end{itemize}
\end{definition}

\begin{lemma}\label{lm:rhoci-rhoimmersion}
Let $\climmersion\colon \oX\to \ovoX$ given by \eqref{eq:oXtooXprime}
be a \rhoimmersion-\rhoci-closed immersion in $\precatAOuB$.
Then \[\rhoSetprime = r(\rhoSet)\amalg e,\]
where $r$ and $e$ are form \Cref{def:precatAOuB:rhoci-closedimmersion,def:rhoirr_rhoimmersion}.
\end{lemma}
\begin{proof}
Lct $\rhoElementprime\in\rhoSetprime$, then 
$\climmersion^*(\rhoElementprime)\in \calS(\schX)$ is either $0$ or in $\rhoSet$.
Define \[\rhoElementprime^\prime=\begin{cases}
\emptyset, \quad \climmersion^*(\rhoElementprime)=0\in \calS(\schX),\\
r(\climmersion^*(\rhoElementprime)), \quad \climmersion^*(\rhoElementprime)\rhoSet \end{cases}
.\] 
Since
\[
\#((e\cup {\rhoElementprime}^\prime))=
\codim_{\ovschX}((e\cup {\rhoElementprime}^\prime))=
\codim_{\ovschX}(\{\rhoElementprime\}\cup (e\cup {\rhoElementprime}^\prime))=
\#(\{\rhoElementprime\}\cup (e\cup {\rhoElementprime}^\prime)),
\]
$\rhoElementprime\in (e\cup \rhoElementprime^\prime)$.
\end{proof}

\begin{remark}\label{def:catAOB:rhoci-closedimmersion}
\emph{\rhoci-closed immersions in $\precatAOB$ and in $\precatAOB$}
are defined by \Cref{def:precatAOuB:retractclosedimmersion} and 
Definitions 
\ref{def:morphisms:precatAOuB->precatAOB} and
\ref{def:morphisms:precatAOB->catAOB}
.
%
\end{remark}


\begin{definition}\label{def:precatAOuB:retractclosedimmersion}
A morphism $\climmersion\colon \oX\to\ovoX$ in $\precatAOuB$ is a 
\emph{retract-closed immersion}
if
$\FschX(\climmersion)$ from \eqref{eq:funct:precatOBtoSch:oX->schX} is 
a retract-closed immersion in $\Sch$, 
\Cref{def:copresheavesSch:retractclosedimmersion}.
\end{definition}

\begin{remark}\label{def:catAOB:retractclosedimmersion}
\Cref{def:precatAOuB:retractclosedimmersion} and \Cref{def:morphisms:precatAOuB->precatAOB}.
define
\emph{retract-closed immersions in $\precatAOB$}.
Then \Cref{def:morphisms:precatAOB->catAOB} defines
\emph{retract closed immersions in $\precatAOB$}.
\end{remark}

\subsubsection{Blow-ups}

\begin{definition}[Blow-up in $\precatAOuB$]\label{def:precatAOuB:BlZoX}
Given $\oX\in \precatAOuB$ 
\eqref{eq:oX} 
and a closed subscheme $Z$ of the scheme $\schX$, 
a blow-up of $\oX$ w.r.t. $Z$ is 
an \preimagegrsup-morphism
\begin{equation}\label{eq:BlZoXinprecatAOuB}
Bl_Z(\oX)\to \oX
\end{equation}
that image along the functor \eqref{eq:funct:precatOBtoSch:oX->schX}
is isomorphic to the morphism $\Bl_Z(\schX)\to \schX$
in $\Sch$.
\end{definition}

\begin{definition}[Blow-up in $\precatAOB$]\label{def:precatAOB:BlZoX}
Given $\ox\in \precatAOB$ 
\eqref{eq:oXtooXprime}, 
\begin{itemize}
\item[(1)]
A blow-up of $\ox$ w.r.t. a closed subscheme $Z$ in $\schX$ 
and 
\item[(2)]
a blow-up of $\ox$ w.r.t. a closed subscheme $\Zprime$ in $\schX^\prime$
\end{itemize}
are 
morphisms 
defined by the commutative squares in $\precatAOuB$
\begin{equation}\label{eq:precatOB:blowups}
\xymatrix{
\Bl_{Z}(\oX)\ar[r]\ar[d]\ar@{}[rd]|{(1)}& \oX\ar[d]\\
\oXprime\ar[r]_{\id_{\oXprime}}& \oXprime
,}
\quad
\xymatrix{
\Bl_{\schX\times_{\schX^\prime}\Zprime}(\oX)\ar[r]\ar[d]\ar@{}[rd]|{(2)}& \oX\ar[d]\\
\Bl_{\Zprime}(\oXprime)\ar[r]& \oXprime
,}
\end{equation}
where the vertical arrows define objects in $\precatAOB$ and the horizontal arrows define morphisms.
\end{definition}


\begin{definition}[\rhoblowup]\label{def:Blrho}
An \emph{\elementaryrhoblowup-morphism}
in $\precatAOuB$
is 
a morphism 
equivalent to 
a bow-up morphism
\eqref{eq:BlZoXinprecatAOuB}
for 
some
$\oX\in \precatAOuB$ given by \eqref{eq:oX}, and
\[Z=Z(\rhoSubset), \quad \rhoSubset\subset\rhoSet.\] 
%
\par\noindent
A \emph{\rhoblowup-morphism}
in $\precatAOuB$
is a composite of a sequence of \elementaryrhoblowup-morphisms.

%

\end{definition}

\begin{definition}\label{def:Monrho_Blrhoelements}
Given $\oX\in \precatAOuB$ \eqref{eq:oX},
denote by $\rhoF$
the set of subsets $\rhoMonSubset$ of $\Monrho$ 
such that $\Bl_{\rhoMonSubset}(\schX)\to \schX$
is a $\rho$-blow-up.

\end{definition}

\begin{definition}\label{def:goodblowup}
A closed subscheme $C$ of a scheme $X$ is {\good } w.r.t a closed subset $Z$ of $X$,
if either $C\subset Z$ or $\codim_{X}(C\cap Z)\geq \codim_X C+\codim_X Z$.

A morphism 
of schemes $\widetilde{X}\to X$
is a {\good } blow-up  w.r.t. a closed subset $Z$ of $X$,
if 
it is isomorphic to 
the blow-up $\Bl_C(X)\to X$ 
for closed subscheme $Z$ of $X$ \good w.r.t. $C$.
\end{definition}

\begin{definition}[\betablowup]\label{def:Blbeta}
An \emph{\elementarybetablowup-morphism}
in $\precatAOuB$
is 
a morphism 
equivalent to 
a bow-up morphism
\eqref{eq:BlZoXinprecatAOuB}
for 
some
$\oX\in \precatAOuB$ given by \eqref{eq:oX}, and
\[Z=Z(\SubsetMonrho), \quad \SubsetMonrho\in\presetsSchrhoBlcup\]
that is a {\good } blow-up w.r.t. the vanishing locus of each element of $\presetsSchrhoBlcup$, \Cref{def:goodblowup},
and 
a \preimagerho-morphism, \Cref{def:rhoirr_rhoimmersion}.
%
%
%
%
\par\noindent
A \emph{\betablowup-morphism}
in $\precatAOuB$
is a composite of a sequence of \elementarybetablowup-morphisms.

%

\end{definition}

\begin{remark}\label{def:catAOB:Blrho}
\emph{\rhoblowup-morphisms} and \emph{\betablowup-morphisms}
in $\precatAOB$ and $\Spc^{\precatAOB}$ 
are defined by Definitions 
\ref{def:morphisms:precatAOuB->precatAOB} and
\ref{def:morphisms:precatAOB->catAOB} 
w.r.t \Cref{def:Blrho} and \Cref{def:Blbeta}.
\end{remark}

\begin{definition}[\FD-object]\label{def:betaFobject}
An \FD-object in $\precatAOuB$ is an object $\oX$ \eqref{eq:oX}
such that
$\presetsSchrhoBlcup\subset \rhoF$.
\end{definition}
\begin{lemma}\label{lm:FDBlbeta}
A {\betablowup } of an \FD-object is an \FD-object.
\end{lemma}
\begin{proof}
$\bl^*(\betaElement) = \bl^*_{\mathrm{proper}}(\betaElement)\cdot \exceptional^\powerexceptional$
\end{proof}

\begin{lemma}\label{lm:anyBlatrhoMonElementistrivialonrhobowup}
Given $\oX\in \precatAOuB$,
for any finite set of finite subsets $\rhoSubsetsCentres$ in $\Monrho$ 
and any set of finite subsets $\rhoSubsetswrtGood$ in $\Monrho$,
there is a finite sequence of elementary 
\rhoblowups
\begin{equation}\label{eq:seqelrhobl}\oX=\oX_0\leftarrow\dots=\oX_{j}\leftarrow\Bl_{Z(\beta_j^\mathrm{el})}(\oX_j)=\oX_{j+1}\leftarrow\dots 
=\oX_l\end{equation}
such that 
each $Z(\beta_j^\mathrm{el})$ is good w.r.t 
inverse images of $Z(\rhoSubsetGood)$ for all $\rhoSubsetGood\in \rhoSubsetsGood$
and  
for each $\rhoSubsetCentre\in \rhoSubsetsCentres$,
the blow-up 
$\Bl_{Z(\bl^*(\rhoSubsetCentre))}(\oX_l)\to\oX_l$
admits a section, 
where $\bl$ is the composite of \eqref{eq:seqelrhobl}.
%
\end{lemma}\begin{proof}
The proof is by a cascade induction with respect to parameters
\begin{itemize}
\item[]
(P1) 
$\#\rhoSubsetsCentres$,
\item[] 
(P2) 
$(
\sum_{\rhoSubsetCentre\in \rhoSubsetsCentres}\#Z(\rhoSubsetCentre), 
\sum_{\rhoSubsetCentre\in \rhoSubsetsCentres}\dim Z(\rhoSubsetCentre)
)$,
\item[] 
(P3) 
$
\sum_{\rhoSubsetCentre\in \rhoSubsetsCentres} 
\sum_{\rhoMonElement\in \rhoSubsetCentre} \deg \rhoMonElement
$%
.
\end{itemize}

(P1 Base)
Suppose $\#\rhoSubsetsCentres=0$. Then the claim holds.

(P1 Step)
Suppose the claim holds for any $\#\rhoSubsetsCentres$ less then the given one.

(P2 Base) 
Suppose 
$Z(\rhoSubsetCentre)=\emptyset$ for each $\rhoSubsetCentre\in \rhoSubsetsCentres$.
Then the claim holds.

(P2 Step)
Suppose 
the claim holds
if 
(P2) is
less than for the given $\rhoSubsetsCentres$.

(P3 Base)
Suppose $\rhoSubsetsCentres$ is a set of subsets of $\rhoSet$.
Let $\rhoSubsetCentreInductionElement\in \rhoSubsetsCentres$.

By induction assumption on (P2)
there is a sequences of elementary $\rho$-blow-ups 
that trivialise
blow-ups with centres of the forms
\[Z(\rhoSubsetCentreInductionElement\cup \rhoSubsetGood),\quad
\rhoSubsetCentre\in \rhoSubsetsCentres, \rhoSubsetGood\in \rhoSubsetsGood,\]
and 
\[Z(\rhoSubsetCentreInductionElement\cup \rhoSubsetCentre),\quad
\rhoSubsetCentre\in \rhoSubsetsCentres\]
that are not equal to $Z(\rhoSubsetCentreInductionElement)$,
and
are \good w.r.t. 
each element of $\rhoSubsetsCentres$ and $\rhoSubsetsGood$.
Composing with 
the blow-up with respect to $\rhoSubsetCentreInductionElement$
reduces the claim to the claim for a smaller (P1). 

(P3 Step)
Suppose the claim would hold 
if 
(P3)
would be less 
than the given one.
Let $\rhoSubsetCentreInductionElement\in \rhoSubsetsCentres$.
\[\rhoSubsetsCentres=(\rhoSubsetCentreInductionElement,\rhoSubsetsCentreshat),\]
and
suppose 
\[\rhoSubsetCentreInductionElement = (\rhoSubsetCentreInductionElementhat,\multa\cdot \multb),\quad a,b\in \Monrho-\{1\}.\]  
By the induction assumption on (P3)
there is a sequence of elementary $\rho$-blow-ups $\seqbl$
that trivialises
blow-up 
with the centre 
$Z(\rhoSubsetCentreInductionElementhat,\multa)$
%
and 
is {\good } w.r.t.  
$\rhoSubsetsGood$
and 
$(\rhoSubsetCentreInductionElementhat,\multa)$.
Then
the claim reduces to the claim for the data
\[\seqbl^*(\rhoSubsetsCentres,\rhoSubsetsGood)=(\seqbl^*(\rhoSubsetCentreInductionElement),\seqbl^*(\rhoSubsetsCentreshat)),\seqbl^*(\rhoSubsetsGood))\]
which is equivalent to the one with 
\[\seqbl^*(\rhoSubsetCentreInductionElement)=
(\seqbl^*(\rhoSubsetCentreInductionElementhat),\seqbl^*(\multa)\seqbl^*(\multb))\]
substituted by
\begin{equation}\label{eq:bproperrhoSetsab}
(\seqbl^*_{\mathrm{proper}}(\rhoSubsetCentreInductionElementhat),\seqbl^*_{\mathrm{proper}}(\multa)\seqbl^*(\multb)),\end{equation}
because $\seqbl$ is {\good } w.r.t. $(\rhoSubsetCentreInductionElementhat,\multa)$.
The claim for \eqref{eq:bproperrhoSetsab} is equivalent to the one
for
\begin{equation}\label{eq:bproperrhoSetsb}
(\seqbl^*_{\mathrm{proper}}(\rhoSubsetCentreInductionElementhat),\seqbl^*(\multb)),\end{equation}
because
$Z(\seqbl^*_{\mathrm{proper}}(\rhoSubsetCentreInductionElementhat),\seqbl^*_{\mathrm{proper}}(\multa))=\emptyset$,
since
$b$ trivialises 
blow-up 
with the centre 
$Z(\rhoSubsetCentreInductionElementhat,\multa)$.
%
%
The claim for 
\eqref{eq:bproperrhoSetsb} holds by the 
induction assumption on (P3). 
\end{proof}

\begin{lemma}\label{lm:anyBlatsetsSchrhoBlElementistrivialonrhoBlsetsofmonomialsrhoElprhoEl}
Given $\oX\in \precatAOuB$,
for any finite subset $\presetsSchrhoBlf$ in $\presetsSchrhoBl$ let
\begin{equation}\label{eq:setsSchrhoBllimElement}
\setsSchrhoBllimElement = \prod_{\presetsSchrhoBlfcupElement\in \presetsSchrhoBlfcup} \presetsSchrhoBlfcupElement
\end{equation}
there is a finite sequence of elementary 
\betablowups
\[\oX=\oX_0\leftarrow\dots=\oX_{j}\leftarrow\Bl_{Z(\beta_j^\mathrm{el})}(\oX_j)=\oX_{j+1}\leftarrow\dots 
=\oX_l\]
such that 
the morphism $\schX_l\leftarrow \Bl_{Z(\setsSchrhoBlElement)}(\schX)\times_{\schX}\schX_l$
adimts a section.
\end{lemma}\begin{proof}
Denote by $\presetsSchrhoBlcup_n$ the set of subsets of $\rhoSet$
formed by unions of $\#\presetsSchrhoBl-n$ different elements in $\presetsSchrhoBl$  
\[\presetsSchrhoBlcup_n=
\{\bigcup \presetsSchrhoBlSubset|\# \presetsSchrhoBlSubset=\presetsSchrhoBl-n\}.
\]
Let \[\setsSchrhoBlElement_n = \prod_{\presetsSchrhoBlfcupElement\in \presetsSchrhoBlcup_n}\presetsSchrhoBlfcupElement\]
Then inverse image of the ideal 
of $\setsSchrhoBlElement_n$
along the blow-up
\[\Bl_{\setsSchrhoBlElement_{n-1}}(\schX)\to\schX\]
equals the product 
of a product of exceptional divisors ideals
and
the product of the vanishing ideals of the inverse images of $Z_{\schX}(\presetsSchrhoBlSubset)$
for all $\presetsSchrhoBlSubset$ such that $\# \presetsSchrhoBlSubset=\presetsSchrhoBl-n$
and the latter inverse images have empty intersection in $\Bl_{\setsSchrhoBlElement_{n-1}}(\schX)$.
So 
the blow-up of $\Bl_{\setsSchrhoBlElement_{n-1}}(\schX)$ 
defined by 
the inverse image of $\setsSchrhoBlElement_n$ is a \betablowup.
\end{proof}

\subsubsection{\'Etale morphisms and infinitesimal neighbourhoods} 

\begin{definition}[\'Etale in $\precatAOB$]\label{def:precatAOuB:Etale} 
A morphism $\etalemorphism$ in $\precatAOuB$ is called \emph{\'etale} if 
$\FschX(\etalemorphism)$ is \'etale, where $\FschX$ is from \eqref{eq:funct:precatOBtoSch:oX->schX}.
The class of \'etale morphisms is denoted $\Et$.
\end{definition}
\begin{remark}[\'Etale in $\catAOB$]
\label{def:catAOB:Etale}
Then the classes $\Et$ called 
\'etale morphisms in $\precatAOB$ and $\catAOB$ are defined by
Definitions 
\ref{def:morphisms:precatAOuB->precatAOB} and
\ref{def:morphisms:precatAOB->catAOB}
w.r.t \Cref{def:precatAOuB:Etale}.
\end{remark}

%
%
%


\begin{definition}[Infinitesimal \'etale neighbourhoods]\label{def:precatAOB:InfinitesimalEtEqNeighbourhood}
The class of morphisms \[\EtInfNeigh=\uLLP(\Et),\] 
Notation \ref{notation:uLLP_AND_uRLP}, is called \emph{infinitesimal \'etale neighbourhoods in $\catAOB$}.
\end{definition}

\begin{definition}[\'Etale in $\Spc(\catAOB)$]\label{def:Spc(catAOB):EtEq}
The class of \emph{\'etale morphisms in $\Spc(\catAOB)$} 
is the class \[\uRLP(\EtInfNeigh),\] Notation \ref{notation:uLLP_AND_uRLP}.
\end{definition}

\subsubsection{\CapitalEtEquational-morphisms}

\begin{definition}[$\EtEq$-morphisms in $\precatAOuB$]\label{def:precatAOuB:EtaleEq} 
A morphism in $\precatAOuB$ is called {\emph{\EtEquational }} 
if the image along the functor 
$\precatAOuBtoSch\colon\precatAOuB\to\Sch$ from \eqref{eq:funct:precatOBtoSch:oX->schX} is 
{\byequations } and \'etale, \Cref{def:Sch:byequations}. 
\end{definition}

\begin{definition}[$\EtEqIRRho$-morphisms in $\precatAOuB$]\label{def:precatAOuB:EtEqIRRho} 
$\EtEqIRRho$-morphisms
in $\precatAOuB$
are
{\etalebyequations }
\rhotriv-morphisms, \Cref{def:rhoirr_rhoimmersion}.
\end{definition}

\begin{remark}[$\EtEqIRRho$-morphisms in $\catAOB$]\label{def:catAOB:EtEqIRRho}
The class of morphisms $\EtEqIRRho$ in $\catAOB$ 
is defined by
Definitions 
\ref{def:morphisms:precatAOuB->precatAOB} and
\ref{def:morphisms:precatAOB->catAOB} 
with respect to Definition \ref{def:precatAOuB:EtEqIRRho}.
\end{remark}

\subsubsection{Morphisms $\ClassImmersionsHBCforLift$}

\begin{definition}\label{def:ClassImmersionsHBCforLift}

$\ClassImmersionsHBCforLift$ denotes 
the class of morphisms in $\precatAOB$ 
that is the image along the functor
\[\precatAOuB\to \precatAOB; \oX\mapsto \id_\oX,\]
of the class of 
\rhoci-closed immersions in $\precatAOuB$ 
of the form \[\climmersion\colon \oX\to \ovoX\]
such that 
$U(\climmersion)$ 
is a closed immersion of essentially smooth local henselian $B$-schemes.


\end{definition}

%


\subsection{$\infty$-category $\catOB$} 
\label{sect:catOB}



\subsubsection{Regularity} 

%

\begin{definition}\label{def:regularintersectionsrhoprecatAOuB}
$\regularstratarhoprecatAOuB$ is the subcategory of $\precatAOuB$ 
spanned by objects of the form \eqref{eq:oX}
such that 
strata of the filtration 
formed by schemes 
$Z(\subsetrhoSet)$ for all subsets $\subsetrhoSet\subset\rhoSet$
are regular.
\end{definition}

\begin{definition}\label{def:regularintersectionsrhoprecatAOB}
$\regularstratarhoprecatAOB$ 
and
$\regularstratarhocatAOB$
are 
the subcategories of 
$\precatAOB$ and $\catAOB$ 
formed by the objects
defined by diagrams
in $\precatAOuB$
that belong to 
$\regularstratarhoprecatAOuB$ of the respective type, \Cref{rem:catAOB}.
\end{definition}



\subsubsection{\rhoPic-objects}

\begin{definition}\label{def:catAOuB:rho_Pic}
An object $\oX\in \precatAOuB$
is a \emph{\rhoPic-object} 
if
there is an element in $\Monrho(\oX)$
that is a section of a very ample line bundle on $\schX$.
%
%
\end{definition}

\begin{definition}\label{def:catAOB:rho_Pic}
\emph{\rhoPic-objects} in 
$\precatAOB$
and
$\catAOB$
are such objects that are
defined by diagrams
of \rhoPic-objects in $\precatAOuB$ of the respective type, \Cref{rem:catAOB}.
\end{definition}


\subsubsection{{\betaLiftcalSringobject}s}
\label{def:betaprodP-object}

\begin{definition}\label{def:betaLiftcalSring-object}
An object $\oX\in \precatAOuB$ is a \emph{\betaLiftcalSringobject} if
for any 
$\presetsSchrhoBlSubset\subset\presetsSchrhoBl$
there is a filtration 
\[
\emptyset = \cupofpresetsSchrhoBlSubsetSubset_0\subset\dots\subset
\cupofpresetsSchrhoBlSubsetSubset_i\subset\dots\subset
\cupofpresetsSchrhoBlSubsetSubset_{n-1}\subset\cupofpresetsSchrhoBlSubset\]
where 
$n=\#\cupofpresetsSchrhoBlSubset$,
such that 
$\#\cupofpresetsSchrhoBlSubsetSubset_i=i$ for each $i=0,\dots ,n-1$,
and
there are retractions
\[r_{i}\colon Z_{\schX}(\cupofpresetsSchrhoBlSubsetSubset_{i-1})\to Z_{\schX}(\cupofpresetsSchrhoBlSubsetSubset_i) ,\quadforall 
i=1,\dots ,n-1.\]
\end{definition}

\begin{definition}\label{def:betaLiftcalSring-object_catAOB}
A \emph{\betaLiftcalSringobject} in $\catAOB$ is an object given by a diagram of {\betaLiftcalSringobject}s in $\precatAOuB$, \Cref{rem:catAOB}.
\end{definition}

\subsubsection{\productprojectivespace}

\begin{definition}\label{def:catAOB:productprojectivespacemorpshism}
A \emph{ {\productprojectivespacemorphism} 
in $\precatAOuB$}
is a morphism $\oX\to \oXprime$ 
that image along the functor $\FschX$ is isomorphic to
\[\schX\simeq \prod_{\alpha\in A}(\PP^{N_\alpha})\times \schXprime\xrightarrow{p} \schXprime\] 
for some set $A$ and integers $N_\alpha$, $\alpha\in A$,
and
$\rhoSet$ is the union of $p^*(\rhoSetprime)$ and 
a subset 
\[\rhoSetP\subset \bigcup_{\alpha\in A} \{ t_{1,i},\dots,t_{N_\alpha,i},t_{\infty,i} \}\subset \calS(\schX) \]
of the set of 
inverse images of the coordinate sections 
$t_{1,i},\dots,t_{N_\alpha,i},t_{\infty,i}\in\calS(\PP^{N_\alpha})$
on the multiplicands $\PP^{N_\alpha}$.
\end{definition}

\begin{definition}\label{def:catAOB:productprojectivespace}
A \emph{ {\productprojectivespace } 
in $\catAOB$} is an object
defined by a diagram, \Cref{rem:catAOB},
in subcategory of the comma-category $\precatAOuB/(B,(\emptyset,,,\emptyset))$
that objects are {\productprojectivespacemorphism}s. %

\end{definition}




\subsubsection{$\infty$-category $\catOB$}

\begin{definition}[$\infty$-category $\catOB$]\label{def:catOB}
The $\infty$-category $\catOB$
is the full subcategory of $\catAOB$
formed by objects of the form
$\limoX$
such that
there is a diagram in 
$\regularstratarhocatAOB$,
\Cref{def:regularintersectionsrhoprecatAOB},
\[\oX\xleftarrow{\rhoblmorphismmu} \BlrhooX\xleftarrow{\rhoblmorphism} \limBlrhooX \xleftarrow{\etalemorphism} \limElimBlrhooX\simeq \limoX\]
such that
\begin{itemize}
\item[(O1)] 
$\oX$
is a
{\rhoPic-\FD
-\betaLiftcalSring-\productprojectivespace }
in $\catAOB$, 
\Cref{def:catAOB:rho_Pic,def:betaFobject,def:catAOB:productprojectivespace,def:betaLiftcalSring-object_catAOB}. 

\item[(O2)]
$\rhoblmorphismmu$ is a composite of a finite set of elementary {\rhoblowup } morphisms,
$\rhoblmorphism$ is a \preimagerho-{\rhoblowup } morphism, \Cref{def:catAOB:Blrho}, 


\item[(O3)]
$\etalemorphism$ is an $\EtEqIRRho$-morphism,
\Cref{def:catAOB:EtEqIRRho}.

\end{itemize}
\end{definition}

\subsection{Localisations and endo-functors}

\subsubsection{$\Blbeta$-localisation and $\Blrho$--localisation} 
\label{sect:Blbeta-localisation}

In this subsection, we define two types of 
localisations
on 
$\Spc(\catOB)$
called
$\Blbeta$-localisation
and
$\Blrho$-localisation.
The definitions and results for the first type, $\Blbeta$, are written below.
The definitions and results 
for $\Blrho$ are the same 
where the role of 
the set $\presetsSchrhoBlcup$
is played 
by 
the set of subsets of $\rhoSet$.

\begin{definition}\label{def:Blbetaequivalence}
An elementary \emph{\betablowup-equivalence} in $\precatAOB$ is 
a \betablowup-morphism given by \eqref{eq:precatOB:blowups}
such that 
\[Z\subset Z_{\schX}(\gr),\quad \Zprime\subset Z_{\schX^\prime}(\gr^\prime).\]
\end{definition}

Let $\Blbeta$ denote the classes of morphisms in $\precatAOB$ and $\catOB$
generated by elementary \betablowup-equivalences and called \betablowup-equivalences.
Let $\Spc_{\Blbeta}(\catOB)$ denote the subcategory of $\Blbeta$-invariant objects in $\Spc(\catOB)$.
Let \[\LbetaBl\colon \Spc(\catOB)\to \Spc_{\Blbeta}(\catOB)\] denote the localisation functor.

\begin{definition}
For any $\oX\in \catOB$,
$\diagramBlrhoeq(\oX)$ is 
the subcategory of 
the comma $\infty$-category 
${\catOB}{/\oX}$ formed by \rhoblowup-equivalences,
and
$\diagramBlbetaeq(\oX)$ is 
the subcategory formed by \rhoblowup-equivalences that are defined by elements of $\setsSchrhoBl$.
\end{definition}
\begin{lemma}\label{lm:Blbeta:cof}
(1) The subcategory of the comma $\infty$-category ${\catOB}{/\oX}$ formed by \rhoblowup-equivalences
is cofinal in
$\diagramBlrhoeq(\oX)$. 
(2) The subcategory of the comma $\infty$-category ${\catOB}{/\oX}$ formed by \betablowup-equivalences
is cofinal in
$\diagramBlbetaeq(\oX)$. 
\end{lemma}
\begin{proof}
The claim for $\beta$-blow-ups follows from \Cref{lm:anyBlatsetsSchrhoBlElementistrivialonrhoBlsetsofmonomialsrhoElprhoEl}.
The claim for $\rho$-blow-ups follows from \Cref{lm:anyBlatrhoMonElementistrivialonrhobowup}.
\end{proof}

\begin{lemma}\label{lm:LBlbeta:FsimeqinjlimF(Blbeta)} 
For any $F\in \Spc(\catOB)$, and any $\oX\in \catOB$,
there is a natural equivalence
\[\LBlbeta(F)(\oX)\simeq \varinjlim_{(\tilde{\oX}\to\oX)\in \diagramBlbetaeq(\oX)} F(\tilde{\oX}).\]
\end{lemma}
\begin{proof}
Because of \Cref{lm:Blbeta:cof},
\[\varinjlim_{(\tilde{\oX}\to\oX)\in \diagramBlbetaeq(\oX)} F(\tilde{\oX})\simeq \varinjlim_{(\tilde{\oX}\to\oX)\in \diagramBlbetaeq(\oX)} F(\tilde{\oX})\]
and
for any \betablowup-equivalence $b^\prime\colon \tilde{\oX}^\prime\to \oX^\prime$ in $\precatAOuB$
and morphism $\oX\to \oX^\prime$
there s a \betablowup-equivalence $b\colon \tilde{\oX}\to \oX$ 
equipped with a commutative square 
\[\xymatrix{
\tilde{\oX}\ar[d]^{b}\ar[r]&\tilde{\oX}^\prime\ar[d]^{b^\prime}\\
\oX\ar[r]&\oX^\prime.}\]
Then
it follows that the class of \betablowup-equivalences in $\catOB$ satisfies similar property.
Then the claim follows.
\end{proof}

The functor
\[\limBlbeta\colon \precatAOuB\to \Spc^{\precatAOuB}\]
that takes $\oX$ to 
\begin{equation}\label{eq:limBlbetaoX}
\limBlbeta(\oX)=\lim \diagramBlbetaeq(\oX)
\end{equation}
induces the functor
\[\limBlbeta\colon \precatAOB\to \Spc^{\precatAOB}\simeq \catAOB\]    
and moreover, the functor
$\limBlbeta\colon \catAOB\to \catAOB$.
The latter functor 
restricts to the functor
\begin{equation}\label{eq:limBlbeta}\limBlbeta\colon \catOB\to \catOB.\end{equation}
\begin{remark}\label{rem:limBlbeta}
By the construction the functor $\Blbeta$ from \eqref{eq:limBlbeta} preserves limits.
\end{remark}
\begin{lemma}\label{lm:limBlbeta:diagramBlbeta}
The arrow $\limBlbeta(\oX)\to \oX$ is the initial object of the $\infty$-category of pro-objects of the $\infty$-category $\diagramBlbetaeq(\oX)$ for any $\oX\in\catOB$.
\end{lemma}
\begin{proof}
The claim holds for any $\oX\in \precatAOB$ by the construction of $\limBlbeta$ from \eqref{eq:limBlbetaoX}. 
Then
it follows for any $\oX\in \catAOB$.
\end{proof}
\begin{lemma}\label{lm:LBlbeta:limBlbeta} 
For any continuous $F\in \Spc(\catOB)$, and any $\oX\in \catOB$,
there is a natural equivalence
\[\LBlbeta(F)(\oX)\simeq F(\limBlbeta(\oX)).\]
\end{lemma}
\begin{proof}
The claim follows from \Cref{lm:LBlbeta:FsimeqinjlimF(Blbeta),lm:limBlbeta:diagramBlbeta}.
\end{proof}

\begin{lemma}
For any $\oX\in \catOB$ such that $\presetsSchrhoBl$ is finite,
$\limBlbeta(\oX))$ is equivalent to $\Bl_{\setsSchrhoBllimElement}(\oX)$ form \eqref{eq:setsSchrhoBllimElement}.
\end{lemma}
\begin{proof}
The claim follow because blow-up on $\Bl_{\setsSchrhoBllimElement}(\oX)$ defined by the inverse image of each element $\setsSchrhoBlElement\in\presetsSchrhoBlcup$ is identity because \eqref{eq:setsSchrhoBllimElement} includes such multiplicant. 
\end{proof}

\begin{remark}\label{rem:limBlrho}
It follows from the construction that the functor $\limBlbeta$ form \eqref{eq:limBlbeta} preserves limits.
\end{remark}



\subsubsection{$(\compintABlbeta)$-localisation} 
\label{sect:compintABlbeta-localisation}

\begin{definition}\label{def:partialDeltaonecompintA}

\[
\begin{array}{lclcl}
\DeltaonecompintA&:=&
\big(
(\PP^1_B, ((t_0,t_\infty), (t_\infty, 1, \{\{t_\infty\}\})) )
\to
(B,\emptyset,,,\emptyset)
\big)
&\in&\precatAOB
,\\
\DeltaonecompintAall&:=&
\mathrm{id}_{
(\PP^1_B, ((t_0,t_\infty), (t_\infty, 1, \emptyset)) )
}
&\in&\precatAOB
,\end{array}\]
where the morphism in the second row is the identity morphism.

\end{definition}
\begin{notation}\label{notation:DeltalcompintA}
$\DeltalcompintA:=(\DeltaonecompintA)^{\times l}=\prod\limits_{\{1,\dots,l\}}(\DeltaonecompintA)$, 
and
$\DeltazerocompintA:=(\DeltaonecompintA)^{\times 0}=(\DeltaonecompintA)^{\times 0}$.
\end{notation}
\begin{definition}\label{def:DeltapcompintA}
$\DeltabonecompintA,\DeltapcompintA\in \catAOuB$ are the objects
\[
\DeltabonecompintA:=
\varprojlim\left(
\xymatrix{\DeltazerocompintA\ar@<0.5ex>[r]\ar@<-0.5ex>[r]&\DeltaonecompintA}
\right).
\]
\[
\DeltapcompintA:=\prod\limits_{\mathbb{Z}_{\geq 1}}(\DeltabonecompintA)\simeq \lim\limits_{l\in\mathbb Z_{\geq 1}} ((\DeltabonecompintA)^{\times l}).
\] 
\end{definition}
\begin{remark}
The image of $\DeltabonecompintA,\DeltapcompintA$ 
is $\catAOB$ is contained in $\catOB$.
\end{remark}
\begin{definition}\label{def:LcompintABlbeta}
The $\infty$-category $\SpccompintABlbeta(\catOB)$
is the subcategory in $\Spc(\catOB)$ spanned by objects
invariant with respect to 
morphisms of the form $\oX\leftarrow \oX\times \DeltaonecompintA$,
and 
\betablowup-equivalences.
%
\[\LcompintABlbeta\colon \Spc(\catOB)\to \Spc_{\compintA,\Blbeta}(\catOB)\]
is the localisation functor.
An $\compintAoneBlbeta$-equivalence
is a morphism that image along $\LcompintABlbeta$ is an equivalence.
\end{definition}

\begin{lemma}\label{lm:LcompintABlbeta}
For any $F\in \Spc(\catOB)$ and all $\oX\in\catOB$,
there is a natural equivalence
\[
(\LcompintABlbeta F)(\oX)
\simeq
\varinjlim_{\Deltatimesb} (\LBlbeta(F))(\oX\times\DeltalcompintA)
.\]
\end{lemma}
\begin{proof}
$\DeltaonecompintAall$
and 
$\DeltaonecompintA$
are intervals in $\catOB[w_{\Blbeta}^{-1}]$, where $w_{\Blbeta}$ is the class of $\Blbeta$-equivalences.
\end{proof}
\begin{lemmacorollary}\label{lmcor:LcompintABlbeta}
For any continuous $F\in \Spc(\catOB)$ and all $\oX\in\catOB$,
there is a natural equivalence
\[
(\LcompintABlbeta F)(\oX)
\simeq
F(\limBlbeta(\oX\times\DeltapcompintA))
.\]
\end{lemmacorollary}

\subsubsection{Endo-functor $F\mapsto F^\Blras$}

\begin{definition}\label{def:rasarrows}

Let 
$\rasarrows$
denote the class of morphisms 
in $\precatAOuB$ 
that image along the functor $\FschX\colon \precatAOuB\to\Sch$ from \eqref{eq:funct:precatOBtoSch:oX->schX} is an isomorphism.
Denote by $\rasarrows^r$ the class of arrows in $\precatAOB$ defined by 
\Cref{def:morphisms:precatAOuB->precatAOB}.
Denote by the same symbol $\rasarrows^r$ the class of arrows in $\catAOB$ defined by 
\Cref{def:morphisms:precatAOB->catAOB}.
Denote by the same symbol $\rasarrows$ the class of arrows in $\catOB$ that belong to the class $\rasarrows$ in $\catAOB$.

\end{definition}

\begin{lemma}\label{lm:rasarrows:catAOB->catOB}
Let $\oX\in \catOB$.
For any $\rasarrows^r$-morphism
\[\rasarrow\colon \rasoX\to\oX\] in $\catAOB$
there is 
a $\rasarrows^r$-morphism
\[{\rasarrow}^\prime\colon {\rasoX}^\prime\to\oX\] in $\catOB$ 
equipped with a $\rasarrows^r$-morphism $\rasoX^\prime\to \rasoX$
that the composite with $\rasarrow$ is equivalent to $\rasarrow^\prime$.
\end{lemma}

\begin{definition}\label{def:FBlras} 
For any $\oX\in \catOB$ 
we consider 
the category 
$\Diagramrasarrows$ that objects are $\rasarrows^r$-morphisms
of the form 
\[\rasoX\to \oX\]
and morphisms given by commutative triangles of $\rasarrows^r$ in $\catOB$.

For any $F\in \Spc(\catOB)$, 
the assertion
\begin{equation}\label{eq:FBlras}
F^\Bl(\oX) 
:=
\varinjlim_{\rasoX\in \Diagramrasarrows} 
\LrhoBl F(\rasoX)
\end{equation}
defines the object 
$F^\Blras\in\Spc(\catAOB)$. 
Moreover, 
the assertion \eqref{eq:FBlras} defines 
the endo-functor on $\Spc(\catOB)$ 
equipped with the natural transformation form the identity one \[F\mapsto F^\Blras,\quad F\to F^\Blras.\]
\end{definition}


\subsubsection{Endo-functor $F\mapsto F^{\geomrealcompactintAK}$}

\begin{definition}\label{def:geomrealKcompintA_AND_geomrealcompactintAK}
Let $K$ be a pointed \cubicalset, \Notationref{def:cubicalset}.
Then 
$\geomrealKcompintA$ is 
the image of $K$
along the 
left Kan extension of the functor
\[\begin{array}{lcl}
(\Deltableqone)^{\mathbb{Z}_{\geq 0}}&\to& \Spc(\catOB)\\
\prod_{i\in \mathbb Z_{\geq 0}}\Delta^{p_i}&\mapsto& \prod_{i\in \mathbb Z_{\geq 0}}\DeltacompintA{p_i}
\end{array}\]
Similarly for $\geomrealcompactintAK$ and $\DeltacompintAall{p_i}$.
\end{definition}

For any {\cubicalset } $K$, 
there is the endo-functor 
$F\mapsto F^\geomrealcompactintAK$
on $\Spc(\catOB)$
such that
there is a natural equivalence
\begin{equation}\label{eq:FgeomrealcompactintAKoX}(F^\geomrealcompactintAK)(\oX)\simeq \Map^\Spc_{\Spc^\bullet(\catOB)}(\oX\wedge\geomrealcompactintAK),F)
\end{equation}

\subsubsection{Properties of $\LcompintABlbeta$: continuous presheaves and endo-functors} 

\begin{lemma}\label{lm:continuous:LcompinABlbeta}
Let $F\in \Spc(\catOB)$ be continuous, then
$\LcompintABlbeta(F)\in \Spc(\catOB)$ is continuous, \Cref{def:continuous,def:LcompintABlbeta}.
\end{lemma}
\begin{proof}
	Since $F$ is continuous, by \Cref{lm:LBlbeta:limBlbeta}
    \[
    (\LcompintABlbeta F)(\oX)
    \simeq
    \varinjlim_{\Delta^l\in\Deltatimesb} F(\limBlbeta(\oX\times\DeltalcompintA))
    \]
    for all $\oX\in\catOB$.
	By 
	\Cref{rem:limBlrho}
	and since $F$ is continuous,
    \[
    F(\limBlbeta(\oX\times\DeltalcompintA))
    \simeq
    \varinjlim_\alpha F(\limBlbeta(\oX_\alpha\times\DeltalcompintA)).
    \]
    Then since colimits commute with colimits, 
    \[
    \varinjlim_{\Delta^l\in \Delta^\bullet} F(\limBlbeta(\oX\times\DeltalcompintA))
    \simeq
    \varinjlim_\alpha \varinjlim_{\Delta^l\in \Delta^\bullet} F(\limBlbeta(\oX_\alpha\times\DeltapcompintA)).
    \]
\end{proof}
%
\begin{lemma}\label{lm:continuous:()geomrealcompactintAK}
Let $F\in \Spc(\catOB)$ be continuous, then
$F^\geomrealallrascompactintAK\in \Spc(\catOB)$ is continuous, \Cref{def:continuous,def:geomrealKcompintA_AND_geomrealcompactintAK}.
\end{lemma}
\begin{proof}
    By \eqref{eq:FgeomrealcompactintAKoX} and since
    \[(\varprojlim_\alpha \oX_\alpha)\times\geomrealcompactintAK\simeq 
    \varprojlim_\alpha (\oX_\alpha\times\geomrealKcompintA),\]
    it follows that
    \[F^\geomrealcompactintAK(\varprojlim_\alpha \oX_\alpha)\simeq \varinjlim_\alpha F^\geomrealcompactintAK(\oX_\alpha).\]
\end{proof}

\begin{lemma}\label{lm:commute:LBlbetacompintA_AND_(-)Bl}
Let $F\in \Spc(\catOB)$.
Then 
$\LBlbetacompintA(F^\Blras)\simeq \LBlbetacompintA(F)^\Blras$.
\end{lemma}
\begin{proof}
It follows form \Cref{def:partialDeltaonecompintA} that 
for any $\oX\in \precatAOB$, 
the subset \[\rhoSet(\oX\times\DeltalcompintA)\subset\calS(\schX\times(\PP^1)^{\times l})\] equals the inverse image of $\rho(\oX)$ along the projection $\schX\times\DeltalcompintA\to \schX$,
the natural morphism is an equivalence
\[(F^{\DeltalcompintA})^\Blras\xrightarrow{\simeq} (F^\Blras)^{\DeltalcompintA}.\]
So the claim follows because $(-)^\Blras$ preserve colimits in $\Spc(\catOB)$. 
\end{proof}
\begin{lemma}\label{lm:commute:LBlbetacompintA_AND_(-)K}
Let $F\in \Spc(\catOB)$.
Then 
$\LBlbetacompintA(F^\geomrealcompactintAK)\simeq (\LBlbetacompintA(F))^\geomrealcompactintAK$.
\end{lemma}
\begin{proof}
The claim follows because
\[(F^\geomrealcompactintAK)^{\DeltalcompintA}\simeq (F^{\DeltalcompintA})^\geomrealcompactintAK,\]
and $(-)^\geomrealcompactintAK$ preserve colimits in $\Spc(\catOB)$. 
\end{proof}

\subsubsection{$\ABHBC$-equivalences}\label{section:ABHBC}
\begin{definition}\label{def:ABHBC:uLLP(S,Et)}
A morphism $\compinttimmersion$ in $\catOB$
is 
an $\ABHBCprop$-equivalence
if 
\begin{itemize}
\item 
$\compinttimmersion\in\uLLP(\calS\to \calL,\calL\to *)$, 
where $*$ is the terminal object of $\Spc(\catOB)$ and $\calS\to\calL$ is from \Cref{def:LIStoLtoT}, 
and 
\item 
$\compinttimmersion$ is an {\infinitesimaletaleneighbourhood } in $\catAOB$, \Cref{def:precatAOB:InfinitesimalEtEqNeighbourhood}.
\end{itemize}
\end{definition}

\subsection{$\trivcofHBC$-invariant objects}
%
%
We recollect some 
results 
and consequences 
on $\ABHBC$-invariant objects from 
\cite{ABHBCinvariance,ABHBClocalisation}
used in the next section.


\begin{theorem}[\protect{\cite{ABHBCinvariance}}]\label{th:ABHBCinvariance:LcompintArhoBl(sinLcapZsequalsZeqmpty)}
For any $F\in \sinLcapZsequalsZeqmpty$, 
$\LcompintABlbeta(F)$ is
$\ABHBC$-invariant.
\end{theorem}

\begin{corollarytheorem}[\protect{\cite{ABHBCinvariance}}]\label{corollarytheorem:Implication:invarianceABHBC->RLPwrtHBC}
Let $F$ be 
continuous and
$\trivcofHBC$-invariant.
Then $\pi_0(F)$ \RLPwrt 
$\ClassImmersionsHBCforLift$-morphisms, \Cref{def:ClassImmersionsHBCforLift}.
\end{corollarytheorem}


%
%
\begin{theorem}[\protect{\cite{ABHBClocalisation}}]\label{corollary:FinvatiantIMPLIESFBlinvariant} 
Let $F\in \Spc(\catOB)$ be invariant w.r.t. $\trivcofHBC$-equivalences. 
Then $F^\Bl$ is invariant w.r.t. $\trivcofHBC$-equivalences. 
\end{theorem}
\begin{theorem}[\protect{\cite{ABHBClocalisation}}]\label{corollary:FinvatiantIMPLIESFKinvariant}
Let $F\in \Spc(\catOB)$ be invariant w.r.t. $\trivcofHBC$-equivalences. 
Then $F^{|S^1|_{\compintAall}}$ is invariant w.r.t. $\trivcofHBC$-equivalences.
where $S^1=\Delta^1/\partial\Delta^1$.
\end{theorem}

\section{Presentation of morphisms in $\HpointedAtopology(\Sm_B)$. }
\label{sect:Presentation}

Throughout the section 
$B\in \Sch$, 
and
$\tau=\Zar$ or $\tau=\Nis$ is the Zariski or Nisnevich topology on $\Sch_B$ or $\Sm_B$.

\subsection{Notation $\HKonSpttau$} 
The fully faithful functor $\Sm_B\to \Sch_B$
induces the fully faithful functors
\[
\FunctorSpcSmBtoSpcSchB\colon
\Spc(\Sm_B)\to \Spc(\Sch_B),\quad
\FunctorHpointedAtopologySmBtoHpointedAtopologySchB\colon 
\SpcAtopology(\Sm_B)\to \SpcAtopology(\Sch_B).
\]
\begin{definition-notation}\label{notation:HKonSpaces}
Given 
$\motSpc\in \SpcAtopology(\Sch_B)$ 
and pointed {\cubicalset } $K$,
\begin{equation}\label{eq:HKonSpaces}
\HKonSpttau(\motSpc)\colon 
\Spc(\Sch_B)\to \Spt; 
\SpcpresheafonSchB\mapsto \Map^{\Spc}_{\SpcAtopology(\Sch_B)}(\SpcpresheafonSchB\wedge K,\motSpc)
\end{equation}
\end{definition-notation}
\begin{notation}\label{notationHonSpaces}
$\HonSpc_\tau(\motSpc):=\HKonSpctau(\motSpc)$ for $K=\Delta^0_+$.
So
$\HonSpc_\tau(\motSpc)(\SpcpresheafonSchB)\simeq 
\Map^\Spc_{\SpcAtopology(\Sch_B)}(\SpcpresheafonSchB,\motSpc)
$.
\end{notation}
\begin{notation}\label{notation:HKonSpaces}
We 
skip
$\motSpc$ or $\tau$ 
for 
$\HKonSpctau(\motSpc)$ from \eqref{eq:HKonSpaces},
when  the context defines 
any 
of them.
\end{notation}
\begin{notation}
We use notation $\HKonSpc$ for 
the presheaf on $\precatAOuB$ and $\precatAOB$ 
that is the inverse image long the functor
$\precatAOB\to \precatAOuB\to \Spc(\Sch_B)$, 
and for the associated continuous presheaf on $\catAOB$,
and its restriction on $\catOB$.
\end{notation}
%
%
%

\subsection{Intersections of line bundle sections and blow-ups}\label{section:Intersections_of_line_bundle_sections_and_blow-ups}

\begin{definition}\label{def:sLonXKZssubsetcupZs}
%
The 
category $\sLonXKZssubsetcupZs$ is a full subcategory in 
$\Set(\catO_B)$ 
that objects are presheaves of sets $F$ on the $\infty$-category 
$\catO_B$
defined by presheaves on $\precatAOuB$ of the form 
\begin{equation}\label{eq:sLonXcapZssubsetcupZs} 
F(\oX) = 
\{
s\in \calS^{\times\powerSoverL}_\calL
| 
\bigcap\limits_{\indexSubsetssrho\in \SetofindexesSubsetsrho_\indexSubsets}Z( s_{\indexSubsetssrho} )\subset 
Z(\sequation_{\indexSubsets}),
\bigcap\limits_{\indexSubsetssrho\in \SetofindexesSubsetsrho_\indexSubsetsempty}
  Z( s_{\indexSubsetssrho} )=\emptyset
\}
\end{equation}
for 
a morphism of finite sets $\powerSoverL\colon \powerS\to\powerL$,
a set $\SetofindexSubSets$
with 
subsets 
and 
elements 
\[
\SetofindexesSubsetsrho_\indexSubsets\subset \powerS,\quad
\sequation_{\indexSubsets} \in (s\cup \Monrho)\]
for each $\indexSubsets\in \SetofindexSubSets$.
\end{definition}


\begin{theorem}\label{th:equiv:sinLonBlcXcapZsiseqmpty->sinLonXcapZssubsetcupZs}
For any $F\in \sinLonXcapZssubsetcupZs$, 
there is a morphism in $\Set(\precatAOuB)$
\[F^\prime \to F\]
that is an equivalence on $\regularstratarhoprecatAOuB$
for some 
$F^\prime\in \Set(\precatAOuB)$ 
of the form 
\[F^\prime=(F^{\prime\prime})^\limBlras,\quad F^{\prime\prime}\in \sinLcapZsequalsZdelta,\] 
\Cref{def:sinLcapZsequalsZdelta}.
\end{theorem}
\begin{proof}

Let $\oX\in\precatAOuB$
.

Let $F$ be given by \eqref{eq:sLonXcapZssubsetcupZs}.
Let 
$F^\prime\in \sinLcapZsequalsZdelta$
be the subpresheaf of $\calS^{\times\powerSoverL}_{\calL}$ 
defined by \eqref{eq:sLonXKZsequalsZdelta}
with 
conditions like in \eqref{eq:sLonXcapZssubsetcupZs} 
with the same $\indexSubsets$ and $\indexSubsetssrho_\indexSubsets$
as in \eqref{eq:sLonXKZsequalsZdelta} for $F$.

The element $f$ is defined by 
a set of sections of line bundles $s$ on $\schX$.
Then inverse images of $s$ 
on the limit of blow-ups $\BloX$ of $\oX$ with centres $Z( s_{\indexSubsets} )$
define an element of $F^\prime$.
By \cite{cciblowupsrho}
there is a finite sequence of 
morphisms 
\[\BloX_n\xrightarrow{b_n}\oX_{n}\xrightarrow{r_n}\dots\xrightarrow{b_l}\oX_l\xrightarrow{r_l}\BloX_{l-1}\xrightarrow{b_{l-1}}\dots\xrightarrow{r_1}\BloX_0\xrightarrow{\simeq}\oX_0=\oX\] in $\precatAOuB$
such that 
\begin{itemize}
\item 
each $b_l$ is a \rhoblowup, 
\item 
each $r_l$ is a \rhoimmersion-morphism such that $\FschX(r_n)$ is an isomorphism in $\Sch_B$ 
and 
\item
the composite
decomposes as
\[\BloX_n\to \BloX\to \oX.\]
\end{itemize}
Hence the inverse images of $s$ defines an element in $F^{\prime\prime}(\BloX_n)$.
The latter element defines an element in $F^\prime(\oX)$.
\end{proof}

\begin{lemma}\label{lm:siLcapZsequalsZdelta<-sinLcapZsequalsZempty}
For any $F\in \sinLcapZsequalsZdelta$, 
there is 
a morphism 
in $\Set(\precatAOuB)$
\[F^\prime \to F\]
that induces an equivalence on 
$\regularstratarhoprecatAOuB$
for some 
$F^\prime\in \Set(\precatAOuB)$ such that
$F^\prime\in \sinLcapZsequalsZempty$, \Cref{def:sinLcapZsequalsZempty}.
\end{lemma}
\begin{proof}
%

Any element $f\in F(\oX)$ 
for any $\oX\in \regularstratarhoprecatAOuB$
is defined by 
the set of irreducible multiplicands of the sections defining the element $f$,
and conversely, the set of irreducible multiplicands is defined by $f$.
The latter sets of irreducible multiplicands form elements $f^\prime\in F^\prime_\alpha(\oX)$
for presheaves $F^\prime_\alpha\in \sinLcapZsequalsZempty$
with the required morphism.
\end{proof}

\begin{corollary}\label{cor:surj:sinLonBlcXcapZsiseqmpty->sinLonXcapZssubsetcupZs}
For any $F\in \sinLonXcapZssubsetcupZs$, 
there is a morphism in $\Set(\precatAOuB)$
\[F^\prime \to F\]
that is an equivalence on $\regularstratarhoprecatAOuB$
for some 
$F^\prime\in \Set(\precatAOuB)$ 
such that
\[F^\prime=(F^{\prime\prime})^\limBlras, 
\quad
F^{\prime\prime}\in \sinLcapZsequalsZempty,
\] 
\Cref{def:sinLcapZsequalsZempty}. 
\end{corollary}
\begin{proof}
The claim follows by the combination of 
\Cref{th:equiv:sinLonBlcXcapZsiseqmpty->sinLonXcapZssubsetcupZs,lm:siLcapZsequalsZdelta<-sinLcapZsequalsZempty}.
\end{proof}

\subsection{Presentation}
\label{section:subsectionPresentation}

\begin{definition}\label{def:bigcupnsLonXKZssubsetcupZsDeltaonecompintAtimesnalphaColimits}
Let 
$
(   
  \bigcup_n(  {(\sLonXKZssubsetcupZs)}^{ (\DeltaonecompintA)^{\times n_\alpha} }  )   )^\Colimits_{\Spc(\precatAOuB)}
$
denote the subcategory of presheaves 
generated by objects of the form
\[
(
F^{\prime}
)^{(\DeltaonecompintA)^{\times n_\alpha}}
, 
\quad 
F^{\prime}\in \sLonXKZssubsetcupZs\] 
via colimits.
\end{definition}


\begin{definition}\label{def:grcomplZar}
A $\grcompZar$-covering of $\oX\in \precatAOuB$
is a subset $\covset$ in $\calS(\oX)$
such that there is the equality $\bigcup_{s\in \covset}(\stackrel{\circ}{X}-Z(s))=\stackrel{\circ}{X}$ of open subschemes in $\stackrel{\circ}{X}=\funcXcomplinfX(\oX)\in \Sch$.
\end{definition}

\begin{lemma}\label{lm:Cov}
Any $\funcXcomplinfX^*(\Zar)$-coverings of $X\in \precatAOuB$ 
admits a refinement by a $\grcompZar$-covering.
\end{lemma}

\begin{definition}\label{def:Cov_Covl}
$\Cov(\oX)$ is the diagram of $\grcompZar$-coverings of $\oX\in \precatAOuB$
with morphisms given by refinements.
$\Cov_l(\oX)$ is the set of 
$\grcompZar$-coverings with $\#\covset =l$, \Cref{def:grcomplZar}.
\end{definition}

\begin{definition}\label{def:LZaroneF}
We define the endo-functor $\LZar$ on $\Spt(B)$ such that
\begin{equation}\label{eq:LZaroneF}\LZar^{(1)}(F)\simeq \varinjlim_{e\in \Cov(X)} \check{C}_{e}(F)\end{equation}
\end{definition}


Let $\check{C}_{e}(F)$ 
denote $U\mapsto \varinjlim_l F(\tilde{X}^{\times l}\times_X U)$
and 
$\check{C}^{\intA}_{e}(F)$ be the 
$\Delta_{\intA}^\bullet$-geometric realisation in $\Spc(\Sch_X)$.
$\check{C}_{e,X}(F)=\check{C}_{e}(F)(X)$,
$\check{C}^{\intA}_{e,X}(F)=\check{C}^{\intA}_{e}(F)(X)$.

\begin{lemma}\label{lm:LZaronel}
There is an equivalence \[\LZar^{(1),l}(F)(X)\simeq \varinjlim_{l} \LZar^{(1),l}(F)(X)\]
where
\[\LZar^{(1),l}(F)(X)= \coprodpointed_{e\in \Cov_l(X)} \check{C}_{e,X}(F).\]
\end{lemma}
\begin{proof}
The claim is provided by the definition of $\LZar$ \eqref{eq:LZaroneF},
and because 
(1) the set $\Cov(X)$ is the union of sets $\Cov_l(X)$
(2) for each $l$, 
the subdiagram of $\Cov(X)$ formed by $\Cov_{l}(X)$ is discrete.
\end{proof}

\begin{lemma}\label{lm:LZaroneLintA}
There is a natural equivalence 
$\LZar^{(1),l}\LintA(F)(X)\simeq \coprodpointed_{e\in \Cov_l(X)} \check{C}^{\intA}_{e,X}(F)$
.
\end{lemma}

\begin{lemma}\label{lm:th:LZarLintAonel}
For any $F\in (\sLonXKZssubsetcupZs)^\Colimits_{\Spc(\precatAOuB)}$, 
$\LZar^{(1)}\LintA(F)(X)$ is in 
$
(   
  \bigcup_n(  {(\sLonXKZssubsetcupZs)}^{ (\DeltaonecompintA)^{\times n_\alpha} }  )   )^\Colimits_{\Spc(\precatAOuB)}
$.
\end{lemma}
\begin{proof}
It 
from \Cref{lm:LZaroneLintA} and \Cref{lm:LZaronel}
that
there are natural equivalences 
\[\LZar^{(1)}\LintA(F)(X)\simeq \varinjlim_{l} \LZar^{(1),l}\LintA(F)(X)
\simeq \varinjlim_{l} \coprodpointed_{e\in \Cov_l(X)} \check{C}^{\intA}_{e,X}(F).\]
Then we note that $(X\mapsto \coprodpointed_{e\in \Cov_l(X)} \check{C}^{\intA}_{e,X}(F))\in \sLonXKZssubsetcupZs$.
\end{proof}



\begin{lemma}\label{lm:th:LZarintA}
For any $F\in (\sLonXKZssubsetcupZs)^\Colimits_{\Spc(\precatAOuB)}$, 
$\LZarintA(X)$ is in 
$
(   
  \bigcup_n(  {(\sLonXKZssubsetcupZs)}^{ (\DeltaonecompintA)^{\times n_\alpha} }  )   )^\Colimits_{\Spc(\precatAOuB)}
$.
\end{lemma}

\begin{proof}
Since 
there is an equivalence
\[\LZarintA(F)\simeq \varinjlim_{p} (\LZar^{(1)}\LintA)^{\circ p} (F),\]
and $HS$ is in 
Hence 
the claim 
for the endo-functor $\LZarintA$ reduces to 
the claim for 
the endo-functors $\LZar^{(1)}\LintA$ 
provided by 
\Cref{lm:th:LZarLintAonel}%
.
\end{proof}

\begin{theorem}\label{th:equiv:Colimits(sLZssubsetcupZs)->Hsup}
Let $\topology = \Zar$.
Let $\motSpc=\SpherewedgelB$. 
Then there is a morphism in $\Spc(\precatAOuB)$
\[
F\to {\funcXcomplinfXquotientXcomplinfXunionSup}^*(\HonSpc),
\]
that is an equivalence 
and 
such that $F$ is in
$(   
  \bigcup_n(  {(\sLonXKZssubsetcupZs)}^{ (\DeltaonecompintA)^{\times n_\alpha} }  )   )^\Colimits_{\Spc(\precatAOuB)}
$.

\end{theorem}

\begin{proof}

For any $\oX\in \precatAOuB$, 
$H(\oX) = \LZarintA(HS)((\schX-Z(\gr))//(\schX-Z(\gr)-Z(\supp)))$,
where $HS = \Map_{}(-,\SpherewedgelB)$,
and
$(\schX-Z(\gr))//(\schX-Z(\gr)-Z(\supp))=\schX\amalg_{\schX-Z(\gr)} (\schX-Z(\gr))\times(\Delta^1_{\intA}/\Delta^0_{\intA})$. 
Since $HS$ is in $\sLonXKZssubsetcupZs$,
$\LZarintA(HS)$ is in 
$
(   
  \bigcup_n(  {(\sLonXKZssubsetcupZs)}^{ (\DeltaonecompintA)^{\times n_\alpha} }  )   )^\Colimits_{\Spc(\precatAOuB)}
$
by \Cref{lm:th:LZarintA}.
Then $H$ is in 
$
(   
  \bigcup_n(  {(\sLonXKZssubsetcupZs)}^{ (\DeltaonecompintA)^{\times n_\alpha} }  )   )^\Colimits_{\Spc(\precatAOuB)}
$
because sections of line bundles on $(\schX-Z(\gr))//(\schX-Z(\gr)-Z(\supp))$
are defined by sections of line bundles on $\schX\times\Delta^1_\intA$.
\end{proof}



\begin{theorem}\label{th:equiv:Colimits(sinLcapZsequalsZempty)->Hsup}
Let $\topology = \Zar$.
Let $\motSpc=\SpherewedgelB$
, $K$ be a pointed cubical set.
Then there is a morphism in $\Spc(\precatAOuB)$
\[F\to {\funcXcomplinfXquotientXcomplinfXunionSup}^*(\HonSpt(\motSpt)).\]
that is an equivalence on $\regularstratarhoprecatAOuB$
and 
such that $F$ is 
 generated by objects of the form
\[
((F^{\prime})^\limBlras)^{(\DeltaonecompintA)^{\times n_\alpha}}, 
\quad 
F^{\prime}\in \sinLcapZsequalsZempty
\] 
via colimits.
\end{theorem}
\begin{proof}
The claim follows form \Cref{th:equiv:Colimits(sLZssubsetcupZs)->Hsup}, \Cref{cor:surj:sinLonBlcXcapZsiseqmpty->sinLonXcapZssubsetcupZs}.
%
\end{proof}

\begin{theorem}\label{th:equiv:Colimits(sinLcapZsequalsZempty)->HKsup}
Let $\topology = \Zar$.
Let $\motSpc=\SpherewedgelB$
, $K$ be a pointed cubical set.
Then there is a morphism in $\Spc(\precatAOuB)$
\[F\to {\funcXcomplinfXquotientXcomplinfXunionSup}^*(\HKonSpt(\motSpt)).\]
that is an equivalence on $\regularstratarhoprecatAOuB$
and 
such that $F$ is in
$(   
  \bigcup_n(  {(\sLonXKZssubsetcupZs)}^{ (\DeltaonecompintA)^{\times n_\alpha} }  )   )^\Colimits_{\Spc(\precatAOuB)}
$, in other words $F$
 is generated by objects of the form
\[
(((F^{\prime})^\limBlras)^{(\DeltaonecompintA)^{\times n_\alpha}})^\geomrealallrascompactintAK, 
\quad 
F^{\prime}\in \sinLcapZsequalsZempty
\] 
via colimits.

\end{theorem}
\begin{proof}
The claim follows form \Cref{th:equiv:Colimits(sinLcapZsequalsZempty)->Hsup}, \Cref{eq:HKonSpaces},
because the endo-functor $(-)^{\geomrealallrascompactintAK}$ on $\Spt(\precatAOuB)$ preserves colimits.
%
\end{proof}

\section{Injectivity Theorem}
\label{section:Injectivity}
\subsection{Lifting propery and Injectivity} 
\label{sect:InvarianceLiftingPorpertyandInjectivity}

\begin{theorem}\label{th:RLP:HsupK} 
Let $\topology = \Zar$, and $B\in \Sch$ be a regular scheme of a finite Krull dimension.
Let $\motSpc\in \SpcAtopology(B)$  
belong 
to the subcategory generated by 
$\PP^{\wedge l}_B$ for all $l\in \mathbb Z_{\geq 0}$ 
via colimits and extensions.
Let $n\in \mathbb Z_{\geq 0}$,
and
$\HSnonSpc=\HSnonSpc(\motSpc)$, \Cref{notation:HKonSpaces}.
\par
Then 
\[{\funcXcomplinfXquotientXcomplinfXunionSup}^*(\HSnonSpc)\in \Spc(\catOBregstratasgrsupp)\]
is $(\ABHBCprop,\compintAone,\Blbeta)$-invariant,
and 
the preshaef 
\[\pi_0({\funcXcomplinfXquotientXcomplinfXunionSup}^*(\HSnonSpc))\in\Set(\catOBregstratasgrsupp)\]
\RLPwrt 
$\ClassImmersionsHBCforLift$-morphisms,
\Notationref{notation:LLP_AND_RLP}, 
\Cref{def:ClassImmersionsHBCforLift},
and \Cref{def:regularintersectionsrhoprecatAOB}. 


\end{theorem}

\begin{proof}

By \Cref{th:equiv:Colimits(sinLcapZsequalsZempty)->HKsup,cor:surj:sinLonBlcXcapZsiseqmpty->sinLonXcapZssubsetcupZs} 
there is an equivalence in $\Spc(\catOB)$ 
\begin{equation}\label{eq:morphism:Fprime->HsupK}
F^\prime\to {\funcXcomplinfXquotientXcomplinfXunionSup}^*(\HSnonSpc)
\end{equation}
such that
$F$ is generated by 
objects of the form
\newcommand{\nd}{N}
\[
F^\prime=(((F^{\prime\prime})^\limBlras)^{(\DeltaonecompintA)^{\times \nd}})^\geomrealcompactintASn,\quad 
F^{\prime\prime}\in 
\Sigma^{\infty}(\sinLcapZsequalsZempty)
\] 
via colimits and extensions.
By \Cref{lm:commute:LBlbetacompintA_AND_(-)Bl,lm:commute:LBlbetacompintA_AND_(-)K}
\begin{equation}\label{eq:LcompintABlrhoFprime}
\LcompintABlbeta(F^\prime)\simeq 
(
(
(
\LcompintABlbeta(F^{\prime\prime})
)^\limBlras
)^{(\DeltaonecompintA)^{\times \nd}}
)^\geomrealcompactintASn).
\end{equation}
By \Cref{th:ABHBCinvariance:LcompintArhoBl(sinLcapZsequalsZeqmpty)} 
$\LcompintABlbeta(F^{\prime\prime})$ is $\ABHBCprop$-invariant.
By \Cref{corollary:FinvatiantIMPLIESFBlinvariant} and \Cref{corollary:FinvatiantIMPLIESFKinvariant} 
the right side presheaf of \eqref{eq:LcompintABlrhoFprime} is $\ABHBCprop$-invariant.
Thus $\LcompintABlbeta(F^\prime)$ is $(\ABHBCprop,\compintAone,\Blbeta)$-invariant.
Then
$\pi_0\LcompintABlbeta(F^\prime)$
\RLPwrt $\ClassImmersionsHBCforLift$-morphisms 
by \Cref{corollarytheorem:Implication:invarianceABHBC->RLPwrtHBC}. 
\end{proof}


\begin{theorem}[\protect{\cite{Morel-connectivity,SS,DHKY,zbMATH07612787}}]\label{th:Injectivitypositiverelativecodimansion}
Let $B\in \Sch$ be a regular scheme of a finite Krull dimension.
Let $U$ be a local henselian essentially smooth scheme over $B$, 
and $Z$ be a closed subscheme of positive relative codimension. 
Let $S\in \Hpointed_{\A^1,\Nis}(B)$, and $K$ be a pointed cubical set, Notation \ref{def:cubicalset}. 
Then the morphism of the pointed sets
\[\supextensionmorphismHANisB\colon 
\MapSet_{\Hpointed_{\A^1,\Nis}(B)}(U_+/(U-Z)_+\wedge K, S)\to \MapSet_{\Hpointed_{\A^1,\Nis}(B)}(U_+\wedge K, S)\]
is trivial.
\end{theorem}

\begin{theorem}\label{th:Injectivity}

Let 
$B\in \Sch$ be regular of Krull dimension one.
Let $U$ be an essentially smooth local henselian scheme over 
$B$.
Let $\rho$ denote a regular function on $B$ that vanishing locus equals a closed point in $B$
as well as its inverse image on $U$.

For any non-negative integer $n$,
the morphism of the pointed sets
\begin{equation}\label{eq:morphism:HsupKUZ->HsupKU}
\supextensionmorphismHAtopologyB
\colon 
\MapSet_{\HpointedAtopology(B)}(U_+/(U-Z(\rho))_+\wedge S^n,\motSpt)\to \MapSet_{\HpointedANis(B)}(U_+\wedge S^n,\motSpt)
\end{equation} 
is trivial
for
for
any 
$\motSpc\in \HpointedAtopology(B)$ 
such that 
$\pi_0{\funcXcomplinfXquotientXcomplinfXunionSup}^*(\HSnonSpc(\motSpc))\in\Set(\catOB)$
\RLPwrt $\ClassImmersionsHBCforLift$-morphisms,
and
$\topology = \Zar$.
\end{theorem}

\begin{proof}
The object 
\[
\objectUrho=\id_{(U,(\{\rho\},(1,\rho,\emptyset))))}\in\catOB\] 
in $\regularstratarhoprecatAOB$, \Cref{def:regularintersectionsrhoprecatAOB}, 
goes by applying the natural transformation of the functors 
${\funcXcomplinfX}\to {\funcXcomplinfXquotientXcomplinfXunionSup}$ 
from \eqref{eq:morphism:funcXcomplinfX->funcXcomplinfXquotientXcomplinfXunionSup}
to
the morphism 
\[U/(U-Z)\to U\in\Spc(\Sch_{B}).\]
Then the morphism \eqref{eq:morphism:HsupKUZ->HsupKU}
is equivalent to to the morphisms
\[
\supextensionmorphismtU
\colon 
{\funcXcomplinfXquotientXcomplinfXunionSup}^*(\HSnonSpt)( \objectUrho )\to
{\funcXcomplinfX}^*(\HSnonSptNis)( \objectUrho )
\]
induced by 
the natural transformation
$\funcXcomplinfX\to \funcXcomplinfXquotientXcomplinfXunionSup$ and morphism
$\HSnonSpc\to\HSnonSpcNis$
defined by \eqref{eq:HKonSpaces}.
\par
Let 
$\rhomap\colon U\to \A^N_U$
be a morphism such that $Z(\rhomap^*(t_0))=Z(\rho)$.
Let 
\[\objectUrhotilde=
(\id_{(\A^N_{U},(\{t_0\},(1,t_0,\emptyset)))})^{h}_{\rhomap(U)} 
,\]
where $(-)^h$ stands for the Henselisation defined as the limit of \'etale neighbourhoods.
The morphism 
\[\timmersionUrho \colon \objectUrho\to \objectUrhotilde
\]
is in $\ClassImmersionsHBCforLift$, \Cref{def:ClassImmersionsHBCforLift}.
By 
assumption
$\pi_0{\funcXcomplinfXquotientXcomplinfXunionSup}^*(\HSnonSpc)$
satisfies RLP w.r.t. $\ClassImmersionsHBCforLift$. 
Then 
for any 
$e\in {\funcXcomplinfXquotientXcomplinfXunionSup}^*(\HSnonSpc)(\objectUrho)$,
there is a lift $\tilde{e}$ of $e$, i.e.
there an element 
\begin{equation}\label{eq:equalityWRTtimmersion:tildeANDe}
\tilde{e}\in {\funcXcomplinfXquotientXcomplinfXunionSup}^*(\HSnonSpc)(\objectUrhotilde)
,\quad
\timmersionUrho^*(\tilde{e})=e.
\end{equation}
Since by \Cref{th:Injectivitypositiverelativecodimansion}, 
the morphism
\[
{\funcXcomplinfXquotientXcomplinfXunionSup}^*(\HSnonSpcNis)( \objectUrhotilde )\to
{\funcXcomplinfX}^*(\HSnonSpcNis)( \objectUrhotilde )
\]
is trivial,
for any $\tilde{e}$ as above,
there is the equality
$\supextensionmorphismtUtilde(\tilde{e})=*$,
where
\[
\supextensionmorphismtUtilde
\colon
{\funcXcomplinfXquotientXcomplinfXunionSup}^*(\HSnonSpc)( \objectUrhotilde )\to
{\funcXcomplinfX}^*(\HSnonSpcNis)( \objectUrhotilde )
.\]
Hence by \eqref{eq:equalityWRTtimmersion:tildeANDe}
$\supextensionmorphismtU(e)=*$
for any $e\in {\funcXcomplinfXquotientXcomplinfXunionSup}^*(\HSnonSpc)(\objectUrho)$.
\end{proof}

\subsection{Acyclicity}\label{sect:Acyclicity} 

\Cref{th:Injectivity} and \cite{ColumnsCousinbiCompexSHfrB}
imply the following consequence.
\begin{corollary}\label{cor:CousinComplex}
Let $B$ be a DVR spectrum.
\par
Let $\motSpectrum\in \SHAZar(\Sm_B)$.
Suppose that
\begin{itemize}
\item[(1)]
$\motSpectrum$ satisfies descent with respect to Nisnevich squares.
\item[(2)] 
$\motSpectrum$ belongs to the subcategory of 
$\SHAZar(\Sm_B)$ 
generated by suspension spectra of 
$\SpherewedgelB$
via shifts and colimits.
\end{itemize}
\par
Then for any essentially smooth local henselian $U$ over $B$
the Cousin complex 
\begin{equation}\label{eq:CousinmotSprectum}
\PresheafmotSpectrum^n(U)\to \bigoplus_{u\in U^{(0)}}\PresheafmotSpectrum_{u}^{n}(U)\to\dots\to \bigoplus_{u\in U^{(c)}}\PresheafmotSpectrum_{u}^{c+n}(U)\to\dots,
\end{equation}
where 
\[\PresheafmotSpectrum^{l}_W(U)=\pi_0\Map_{\SH_{\A^1,\Nis}(B)}(U/(U-W), L_{\A^1,\Nis}(\motSpectrum)[l]),\]
is acyclic on $U$.
\end{corollary}
\begin{proof}
By \Cref{th:Injectivity}
the complex 
\[\PresheafmotSpectrum^n(U)\to \PresheafmotSpectrum_{\eta\times_B U}^{n}(U)\to\PresheafmotSpectrum_{z\times_B U}^{1+n}(U)\to\dots,\]
where $z\in B$ is the image of the closed point of $U$, and $\eta\in B$ be the image of the generic point,
is acyclic.
The morphism from 
the Cousin complex \eqref{eq:CousinmotSprectum}
to the complex 
\[\PresheafmotSpectrum^n(U)\to \PresheafmotSpectrum_{\eta\times_B U}^{n}(U)\to\PresheafmotSpectrum_{z\times_B U}^{1+n}(U)\to\dots,\]
to the Cosuin complex 
is 
a quasi-isomorphism
by \cite{ColumnsCousinbiCompexSHfrB}
for $U=X^h_x$ for any $X\in \Sm_B$, $x\in X$.
Consequently, 
the morphism is a quasi-isomorphism for any essentially smooth local Henselian $U$ because
the Cousin complex is compatible with affine transition morphisms of schemes.
Hence the \eqref{eq:CousinmotSprectum} is acyclic.
\end{proof}

Further, using \cite{Morel-Voevodsky,zbMATH01194164,zbMATH06773295}, we 
deduce the claim of the Gersten conjecture for $K$-theory for local henselian essentially smooth schemes $U$ over a one-dimensional regular base scheme $B$. 
\begin{theorem}\label{th:Ktheory}
    Let $B$ be a DVR spectrum.
    Let $U$ be
    an essentially smooth local henselian scheme over $B$.
    The complex
    \eqref{eq:KthComplexUn}
    is acyclic.
\end{theorem}
\begin{proof}

The K-theory presheaf is isomorphic to 
the presheaf of homotopies of $S^1$-spectra presheaf $\mathrm{K}$ on $\Sm_B$ 
associated by Secal's machine 
to the $\infty$-commutative monoid formed 
by nerves of the groupoids of  
vector bundles $\Vect=\bigoplus_{L}\Vect_{L}$,

T
he 
presheaves $\Vect_{L}\in \Spc(\Sm_B)$ 
satisfy Nisnevich descent by \cite[Proof of Th 5.2.3]{zbMATH06773295}.
Hence $\mathrm{K}\in \Spc(\Sm_B)$ satsfies Nisnevich descent.
Since presheaves 
$K_n$ are $\A^1$-invariant on smooth 
affine $B$-schemes for all $n\geq 0$
by \cite[Prop 4.1]{zbMATH05778592}
and Poincar´e Duality \cite[Section 7.1]{zbMATH05778592}%
,
the presheaf $\mathrm{K}$ is $\A^1$-invariant.
Then
it follows form 
by \cite[Th 5.1.3]{zbMATH06773295} 
that $\LAZar(\mathrm{K})\simeq \LANis(\mathrm{K})$ in $\Spt(\Sm_B)$
and 
\begin{equation}\label{eq:KnpinAZarVectPinANisVect}K_n(X)\simeq \pi_n(\LAZar(\mathrm{K}))\simeq \pi_n(\LANis(\mathrm{K})).\end{equation}

Because of Bott periodicity, \cite[Th 6.8]{zbMATH01194164}, and by \eqref{eq:KnpinAZarVectPinANisVect}, 
there is a $\PP^1$-spectum 
\begin{equation}\label{eq:KPSpectrum}(\LAZar(\mathrm{K}),\LAZar(\mathrm{K}),\dots )\in \Spc(\Sm_B)[\PP^{\wedge -1}_B]\end{equation}
that represents K-theory.
Since by \cite[Prop 3.13, 4.14]{Morel-Voevodsky}
\[K_n(X)\simeq \MapSet_{\Hmot_{\A^1,\Nis}}(X\wedge S^n,\Gr_{\infty,\infty}\times B),\] 
the schemes $\Gr_{N,L,B}$ are cellular for $N, L\in \mathbb Z_{\geq 0}$, 
the spectrum \eqref{eq:KPSpectrum} is generated by $\PP^l_B$ for $l\in \mathbb Z$. 

Then the claim follows for 
by \Cref{cor:CousinComplex}
applied to \eqref{eq:KPSpectrum}.
\end{proof}

\section{Construction of the $\phi$-motivic homotopy category}
\label{sect:Motivichomotopycategories}

In this section we present the definition of $\phi$-motivic homotopy category,
and outline its properties.

\subsubsection{Motivic category of ``compactified'' schemes}

\begin{definition}\label{def:SchCompactified}
$\SchCompactified_B$ is the category with objects
\[(\schX,\gr),\]
where
\begin{itemize}
\item
$\schX\in \Sch_B$,
and 
$\gr\in\calS(\schX)$
invertible at each generic point of $\schX$
and admits locally a presentation as a product $\prod_{\alpha}\gr_\alpha$ for some 
$\gr_\alpha\in \calS(\schX)$
such that 
for any subset $A^\prime$ of $A$, the connected components of 
$\cap_{A^\prime\subset A}Z(\gr_\alpha)$ are irreducible and have codimension $\#(A^\prime)$ in $\schX$,
\end{itemize}
and morphisms 
\[(\schX,\gr)\to (\schX^\prime,\gr^\prime)\]
defined by
\begin{itemize}
\item
a morphism $f\colon \schX\to\schX^\prime$ in $\Sch_B$
such that 
\item
$\gr\leq f^*(\gr^\prime)$.
\end{itemize}
\end{definition}
%

%
%

\begin{definition}\label{def:motiiccategorySchcompB}
The $(\topology,(\PP^1,\infty),\Bl)$-motivic homotopy category of compactified schemes over $B$
is 
the $\infty$-category 
\[
\Spcpointed_{\topology,(\PP^1,\infty),\Bl}(\SchCompactified_B)
\]
that is the subcategory of $\Spc(\SchCompactified_B)$, 
\Cref{def:SchCompactified},
spanned by 
$(\topology,(\PP^1,\infty),\Bl)$-invariant objects.
Here
\begin{itemize}
\item
$\topology$ indicates local equivalences with respect to the inverse image of the topology $\topology$ on $\Sch_B$ along the functor from \Cref{def:FunctorSchCompactified};
\item
$(\PP^1,\infty)$ indicates invariance with respect to morphisms \[(\schX,\gr)\leftarrow (\schX,\gr)\times (\PP^1,\infty);\] 
\item
$\Bl$ indicates the invariance with respect to 
the blow-ups of objects $(X,\gr)$ 
that centers $C\subset X$ 
that have locally on $X$ 
a presentation
$C=Z(\eqC_0,\dots\eqC_l)$ for 
some sections $\eqC_\alpha$ of $\calS$ 
such that 
$\codim_X C=l+1$,
and 
$\eqC_0\leq \gr$.
\end{itemize}

\end{definition}

\begin{remark}
The argument of 
\Cref{th:Injectivity} 
applies to 
$\motSpt$ in 
$\Spcpointed_{\Zar,(\PP^1,\infty),\Bl}(\SchCompactified)$ 
generated by $\PP^{\wedge l}_B$ via shifts and colimits
and
any
\[(U,\gr), \quad Z\subset U,\]
so that
$U$ is local henselian 
and
$(U,Z(\prod I(\gr)\cdot I(Z))\in \SchCompactifiedRegular$, 
where $I(\gr)$ and $I(Z)$ denote the sheaf of ideals defined by the divisor $\gr$ and the closed subscheme $Z$,
and $\SchCompactifiedRegular$ is the subcategory of objects $(\schX,\gr)\in \SchCompactified$ such that $\schX$ is regular, and $Z(\gr)$ is regular.
\end{remark}

\subsubsection{$(\Blrho,\compintAone,\ABHBC,\Blrasequiv,\topology)$-motivic homotopy category}
\newcommand{\catOsuppunitB}{\catOB^{\supp=1}}

\begin{definition}\label{def:BlrhocompintAoneABHBCClrasequivtopology-motivichomooptycategory-equivalences}
We define the following classes of morphisms 
in the category $\precatAOuB$, and consequently, $\catOB$
:
\begin{itemize}
\item $\tau$-local equivalences, where $\tau$ is the topology on $\precatAOuB$ induced by a topology $\tau$ on $\Sch_B$
as the inverse image along the functor $\FunctorSch$,
\item
$\arrowequiv $-equivalences defined by the counit of the 
adjunction 
$\precatAOB\simeq (\precatAOuB)^{\Delta^1} \leftrightarrows \precatAOuB$,
i.e. morphisms $
\vec{Y}
\to \ox$ such that the induced morphism $
Y
\to\oX$ is an equivalence,
\item
$\Blrasequiv$-equivalences are defined by 
morphisms $\ox$ given by \eqref{eq:oXtooXprime} in $\precatAOuB$ such that
$\FschX(\ox)$ is an isomorphism and
$\gr = \FschX(\ox)^*(\gr^\prime)$,
\end{itemize}
and define the following class of morphisms 
in the $\infty$-category $\Spc(\catOB)$%
:

\begin{itemize}
\item
$\supp$-equivalences
are morphisms 
$\mathrm{cofiber}\to \oX$, where $\oX\in \precatAOuB$ is given by \eqref{eq:oX}, and $\mathrm{cofibre}=\operatorname{cofib}((\schX,(\rho,(\gr\cdot\supp,1,\presetsSchrhoBl)))\to (\schX,(\rho,(\gr,1,\presetsSchrhoBl))))$.
\end{itemize}

\end{definition}

\begin{definition}\label{def:BlrhocompintAoneABHBCClrasequivtopology-motivichomooptycategory}

$(\Blrho,\compintAone,\ABHBC,\Blrasequiv,\arrowequiv,\topology)$-motivic homotopy category defined as 
the subcategory
\[\Spc_{\topology,\supp,\arrowequiv, \Blrasequiv,\ABHBC,\compintA,\Blrho}(\catOsuppunitB)\]
of 
$(\topology,\arrowequiv, \Blrasequiv,\ABHBC,\compintAone,\Blrho)$-invariant objects
in $\Spc(\catOsuppunitB)$. 
\end{definition}
\begin{definition}\label{def:FunctorSchCompactified}
$\FunctorSchCompactified$ is the functor
\[\begin{array}{lclcl}
\FunctorSchCompactified&\colon& 
\precatOB&\to&\SchCompactified_B\\
&&(\ox\colon \oX\to\oXprime)&\mapsto&
(\schX,\gr)
.\end{array}\]
\end{definition}
\begin{definition}\label{def:funcXcomplinfXquotientXcomplinfXunionSupCompactified}
$\funcXcomplinfXquotientXcomplinfXunionSupCompactified$ is the functor
\begin{equation}\label{eq:funcXcomplinfXquotientXcomplinfXunionSupCompactified}\begin{array}{lclcl}
\funcXcomplinfXquotientXcomplinfXunionSupCompactified&\colon& 
\precatOB&\to&\Spc(\SchCompactified_B)\\
&&(\schX,(\rho,(\gr,\supp,\presetsSchrhoBl)))&\mapsto&
(\schX,\gr)/(\schX,\gr\cdot\supp)
.\end{array}\end{equation}
\end{definition}
\begin{remark}\label{rem:BlrhocompintAoneABHBCBlrasNis-motivichomotopycategory:Injectivity}
The claim of \Cref{th:Injectivity} holds for 
$(\Blrho,\compintAone,\ABHBC,\Blrasequiv,\Nis)$-motivic homotopy category, \Cref{def:BlrhocompintAoneABHBCClrasequivtopology-motivichomooptycategory}, in the sense that
the fibre of 
the map
\[F( (U,(\rho,1,\rho)) )\to F( (U,(\rho,1,1)) )\]
is trivial
for any $F\in \Spc_{\Nis,\Blrasequiv,\ABHBC,\compintA,\Blrho}(\catOB)$. 
\end{remark}

\subsection{Sketch for Comparison Theorem}\label{sect:ComparisonEq} 

%

The subcategory of continuous objects in 
\[\Spc_{\topology,\supp,\arrowequiv, \Blrasequiv,\ABHBC,\compintA,\Blrho}(\catOB)\]
is equivalent to 
the $\infty$-category
\[\Spc_{\tau,(\PP^1,\infty),\Bl}(\Sch^{\mathrm{t}\Compactified}_B).\]
Here
$\Sch^{\mathrm{t}\Compactified}_B$
is the intersection $\SchComp_B$
and the subcategory 
$\limSchComp_B$ of 
$\Spc^{\SchComp_B}$ defined by analogous conditions as for $\catOB$ in \Cref{def:catOB}.
We outline the proof.

Firstly, 
we introduce notation $\catOsuppunitB$ 
for the subcategory of $\catOB$ 
spanned by the objects defined by diagrams in $\precatAOuB$ 
of objects of the form \eqref{eq:oX}  such that 
$\supp=1$ is the trivial element in $\calS(\schX)$,
For any $\tau\supset \Zar$,
the adunction \begin{equation*}\label{eq:catOsuppBleftrightarrowscatOB}\begin{array}{lcl}\catOsuppunitB&\leftrightarrows& \catOB\\
(\schX,(\rhoSet,(\gr,1,\presetsSchrhoBl)))&\mapsfrom& (\schX,(\rhoSet,(\gr,1,\presetsSchrhoBl))),\\
(\schX,(\rhoSet,(\gr,1,\presetsSchrhoBl)))&\mapsfrom& (\schX,(\rhoSet,(\gr,\supp,\presetsSchrhoBl))).\end{array}\end{equation*} 
induces an equivalence
\[
\Spc_{\topology,\supp,\arrowequiv, \Blrasequiv,\ABHBC,\compintA,\Blrho}(\catOB)
\simeq 
\Spc_{\topology,\arrowequiv, \Blrasequiv,\ABHBC,\compintA,\Blrho}(\catOsuppunitB).\]

Further, 
there is an equivalence 
\begin{equation}\label{eq:ComparisonTheorem}
\Spc_{\tau,\arrowequiv,\Blrasequiv,\compintA,\Blrho}(\catOsuppunitB)
\simeq
\Spc_{\tau,(\PP^1,\infty),\Bl}(\limSchComp_B)
\end{equation}
of the localisations of the categories of presheaves on $\catOB$ and $\limSchComp_B$
with respect to the equivalences
of the types mentioned in the subscripts.
Without providing the required definitions, 
we explain that
the arrows category $\precatAOB= \precatAOuB^{\Delta^1} 
$ and the data $\rhoSet$ 
can be reduced
in the construction of the left side 
because $\arrowequiv$-equivalences equalise $\precatAOuB$ and $\precatAOB$,
and because 
$(\arrowequiv,\Blrasequiv)$-equivalences allow 
modifying $\rhoSet$ in arbitrary way.
\printbibliography
\end{document}